\documentclass[12pt,a4paper]{article}
\usepackage[dvips]{graphicx}
\usepackage{setspace}
\usepackage{amsmath}
\usepackage{overcite}
\usepackage{amsbsy}
\usepackage{amssymb}
\usepackage[utf8]{inputenc}
\usepackage[T1]{fontenc}

\newcommand{\ud}{\ensuremath{\mathrm{d}}}

\newcommand{\bmi}{\ensuremath{\boldsymbol}}

\newcommand{\bnabla}{\ensuremath{\boldsymbol\nabla}}
\newcommand{\bcdot}{\ensuremath{\boldsymbol\cdot}}

\newcommand{\p}{\ensuremath{\partial}}
\newcommand{\beq}{\begin{equation}}
\newcommand{\eeq}{\end{equation}}
\newcommand{\eps}{\ensuremath{\epsilon}}

\long\def\symbolfootnote[#1]#2{\begingroup%
\def\thefootnote{\fnsymbol{footnote}}\footnote[#1]{#2}\endgroup}

\begin{document}

\begin{center}
{\bf \Large Effects of both diffuse and oblique collimated irradiation on phototactic bioconvection}\\[12pt]
{M. K. Panda$^1$\symbolfootnote[1]{Corresponding author;
e-mail:mkpanda@iiitdmj.ac.in}, Shubham Kumar Rajput$^2$ }\\
{$^1$ Department of Mathematics, PDPM Indian Institute of Information Technology Design and Manufacturing, Jabalpur 482005, India}%\\
%{$^2$ Department of Mathematics, Central University of Himachal Pradesh, TAB, Shahpur, District Kangra
%(H.P.), India, Pin: 176206.}  
\end{center}

\noindent
{\bf Abstract}

The linear stability of a finite-depth algal suspension is investigated numerically with particular emphasis on the effects of angle of incidence. The suspension of phototactic algae is uniformly illuminated by both diffuse and oblique collimated irradiation. The bioconvective solutions show a transition of the most unstable mode of disturbance from the stationary (overstable) to overstable (stationary) state at the variation in angle of incidence for fixed parameter ranges. Furthermore, a transition from mode 2 to mode 1 instability is noticed for some parameter values as the angle of incidence varies. Oscillatory modes of disturbance are also predicted at the increment in angle of incidence (or cell swimming speed).

%\keywords{impinging at an off-normal angle at the top.}
\newpage

\section{INTRODUCTION}
\label{sec1}

The phenomenon of spontaneously formation of patterns in shallow suspensions of randomly, but on an average upwardly swimming micro-organisms which are a little denser than the medium is named as bioconvection \cite{pk:kp,hp:ph}. The mostly found micro-organisms that participate in bioconvection with the above features are  bacteria and algae. The generated patterns in bioconvection also disappear, when the micro-organisms stop swimming. It is well known that there are examples where up-swimming and higher density are not involved in the process of pattern formation \cite{pk:kp}. Micro-organisms swim in particular directions on an average due to their response to certain stimuli called \textit{taxes}.  The well recognized examples of \textit{taxes} are \textit{gravitaxis, gyrotaxis, phototaxis} etc. The respond to gravity or acceleration is referred to as gravitaxis and the negative \textit{gravitaxis} denotes swimming vertically upwards. \textit{Gyrotaxis} is defined as the balance between the torque due to gravity and viscous forces (arising from local shear) for a bottom-heavy micro-organism. \textit{Phototaxis} denotes the swimming in the direction of light intensity gradient vector, while start/stop swimming behaviour is observed in photophobic response. Self-shading is a mechanism by which the algae absorb the light incident on them via photosynthesis and scatter it \cite{vh:hv}. This article is relevant to phototaxis only.

It is shown from the experiments that different types of illumination  intensity (e.g., diffuse and/or oblique (vertical) collimated irradiation) and their magnitude may significantly affect the  patterns evolved in bioconvection \cite{wa,ka85,ka86,vin95, kh_1997, williams_11, kage_13,Mendelson98,Kitsunezaki2007}. Bright (strong) light destroys steady patterns in suspensions of micro-organisms or prevents formation of patterns in well-stirred cultures. Light intensity can also modify the shape, size, structure, symmetry and/or scale of the pattern \cite{ka85,williams_11,ka89}. The following reasons may be responsible for it. First, the motile phototactic algae obtain energy through their  photosynthetic pigments (e.g., chlorophyll and carotenoid) and their swimming trajectories can be modified via phototaxis. They swim towards the light source (positive phototaxis) when the light intensity $G$ is below a critical value $G_c$, and swim away from the light source (negative phototaxis) when $G>G_c$ \cite{had87}. The algae therefore remain at a location where $G\approx G_c$.  The second reason for pattern change may be due to modification in their collective behavior for the available light and it is explained as follows. It is well recognized that the algae absorb and scatter the light incident on them and the absorption  results the light intensity $G$ to decrease monotonically as the oblique collimated irradiation propagates through the algal suspension. But, diffuse irradiation propagates more uniformly across the suspension than the oblique collimated (direct) one. Thus, self-shading is less effective under diffuse irradiation due to incomplete movement of chloroplasts (eyespots) in their environment unlike the oblique collimated irradiation. The effect of scattering appears to be twofold: it decreases $G$ at a point by deflecting light away from the incidence trajectory whereas, it increases $G$ at that point due to contributions to intensity by scattering from the elemental volume of algal suspension.

\begin{figure}[!h]
\begin{center}
\includegraphics[scale=0.5]{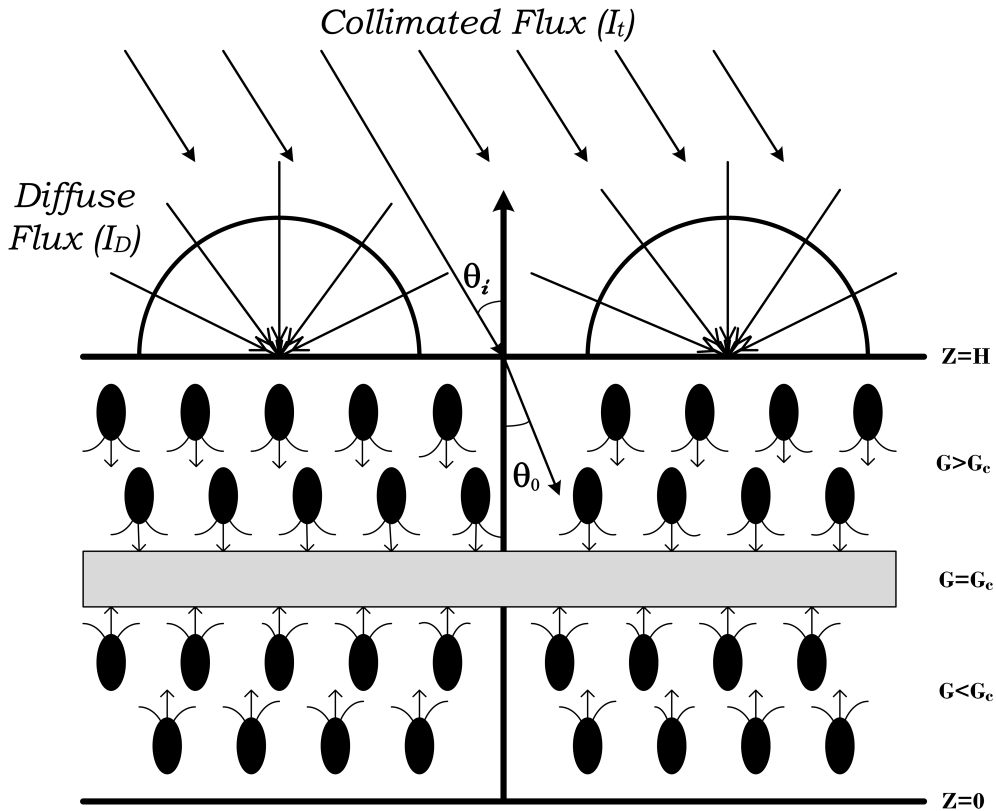}
\end{center}
\caption{Formation of the sublayer in the interior of the suspension. Here $G_c$ is the
critical total intensity.}
\label{int}
\end{figure}

The phototaxis model proposed by Panda \textit{et al.} \cite{prmm16} is used here. In their model, the governing system for bioconvection consists of the Navier--Stokes equations for an incompressible fluid coupled with a conservation equation for algae and the radiative transfer equation (hereafter referred to as RTE) to govern transport of oblique light. Also, the algal suspension  is illuminated by both diffuse and vertical (i.e. not oblique) collimated irradiation. In a natural environment, the sun strikes the surface mostly at different off-normal angles. Thus, the diffuse irradiation is caused by scattering of oblique  collimated irradiation by the effects inside atmosphere (i.e. presence of water droplets (clouds) etc.) \cite{np}. Eventually, the light intensity profiles may be redistributed across the algal suspension via the radiation field and they control the photosynthesis via phototaxis. Furthermore, the realistic estimates of solar radiation transfer across the algal suspension may also affect the time required for formation of patterns in bioconvection \cite{Gittleson1968}. An oblique collimated irradiation may also influence the rate of formation and distribution of aggregations of algae in bioconvection \cite{Gittleson1968}. Solar radiation transfer in algal suspensions relevant to solar energy utilization (particularly in designing efficient photo-bioreactors) and bioconvection due to phototaxis do not appear to be completely understood. Also, the modeling capability for estimating the appropriate radiation field in an algal suspension illuminated by both diffuse and oblique collimated irradiation in a natural environment is not available \cite{gp,williams_13,Bees_14, np}. Thus it is necessary to include oblique collimated irradiation in addition to diffuse irradiation in realistic and reliable models of phototaxis to describe the swimming behavior of algae accurately \cite{kh_1997, macros, incropera_81,incropera_77}. Panda \textit{et al.} \cite{prmm16} calculated the radiation field by neglecting the oblique collimated irradiation in their phototaxis model. In contrast to Panda \textit{et al.} \cite{prmm16}, here the finite-depth algal suspension is illuminated uniformly by both diffuse and oblique collimated irradiation  [see Sec. \ref{sec2}].

Consider a dilute suspension of phototactic algae illuminated by both oblique collimated and diffuse irradiation. To study the bioconvective instability for  such a suspension, the basic state is the one where there is a balance between phototaxis and diffusion of cells. The illuminating irradiation is attenuated (absorbed and scattered) as it travels across the algal suspension.
This results in up swimming (positive phototaxis) in the lower regions ($G < G_c$) and down swimming (negative phototaxis) in the upper regions ($G > G_c$) of the fluid by the algae. Eventually, a horizontal, concentrated layer of algae (the sublayer) is formed at a position where $G=G_c$ in the basic steady state [see Fig.~\ref{int}]. The position of the sublayer is close to the top (bottom) of the algal suspension for low (high) intensities. Since the region below (above) the sublayer is gravitationally unstable (stable), the fluid motions in the unstable layer penetrate the upper stable layer if the fluid layer becomes unstable. This phenomenon occurs in a wide variety of convection problems as an example of penetrative convection \cite{bst}.

Vincent and Hill \cite{vh:hv} investigated bioconvection  in a suspension of phototactic algae. They performed a linear stability analysis of the basic equilibrium solution and found stationary and oscillatory modes of disturbance. Ghorai and Hill \cite{ghp}  simulated numerically two-dimensional phototactic bioconvection using the model proposed by Vincent and Hill \cite{vh:hv}.
 But both of these studies were examined for a non-scattering algal suspension. Ghorai \textit{et al.} \cite{gph} investigated the onset of bioconvection via linear stability theory by assuming that the scattering by phototactic algae is isotropic. A bimodal steady-state profile has been noticed in their study due to isotropic scattering by algae for some parameters. Examples of oscillatory modes of disturbance were also found by them for certain ranges of parameters. Ghorai and Panda \cite{gp} examined the onset of bioconvection via linear stability in an forward scattering suspension of phototactic algae. They observed a transition from a stationary (oscillatory) to an oscillatory (stationary) mode with the variation in the forward scattering coefficient for certain values of the parameters. 
Panda and Ghorai \cite{pg} examined nonlinearly phototactic bioconvection in an absorbing and isotropic scattering suspension confined in a two-dimensional geometry. They legitimated that the obtained patterns differ qualitatively from those found by Ghorai and Hill \cite{ghp} at a higher critical wavelength due to the effects of scattering. Afterwards, Panda and Singh \cite{pr} simulated two-dimensional phototactic bioconvection confined with rigid sidewalls using the model proposed by Vincent and Hill \cite{vh:hv} in $xz$--plane. A significant stabilizing effect on suspension due to lateral rigid walls has been observed in their study for some governing parameters. However, these studies did not include the effects of diffuse irradiation. Panda \textit{et al.}  \cite{prmm16}  investigated the effects of diffuse irradiation on an isotropic scattering algal suspension and the diffuse irradiation had significant stabilizing effect as reported by them. Also, the transition of bimodal vertical concentration profiles at base state to unimodal ones were observed due to diffuse irradiation. Again, the impact of diffuse irradiation was significant on the critical states (i.e., Rayleigh number and wavenumber) at suspension instability in contrast to the collimated irradiation alone. Panda \cite{mkp20} investigated the effects of forward scattering on the bioconvective instability with both diffuse and collimated irradiation. He observed about the transition of bimodal base concentration profiles into the unimodal ones and vice versa for certain parameters due to forward scattering when absorption (self-shading) is insignificant. Furthermore, the most unstable solution shifts from mode $1$ to mode $2$ when the absorption (self-shading) is significant and a single oscillatory branch bifurcates (disappears) from the stationary branch for some parameters. However, the effects of oblique collimated irradiation were not incorporated in the aforesaid studies. First time, Panda \textit{et al.} \cite{PPS} examined the effects of oblique collimated light on bioconvection in a non-scattering algal suspension. The solutions show a transition of the most unstable mode from stationary (overstable) to overstable (stationary) state for certain parameters at the variation in angle of incidence. A transition from mode 2 to mode 1 instability was also observed at the variation of angle of incidence. More recently, Kumar \cite{sandeep} investigated the effects of oblique collimated irradiation on bioconvection in an isotropic scattering algal suspension. He found that the bioconvection solutions are generally oscillatory (or stationary) and more stable (or unstable) in the case of a rigid (or stress-free) upper boundary. However, no study on the onset of phototactic bioconvection that incorporates the effects of both oblique collimated and diffuse irradiation on an algal suspension has been hitherto carried out. Therefore, the effects of oblique collimated irradiation on bioconvection are investigated in the same vicinity.

\begin{figure}[!h]
\begin{center}
\includegraphics[scale=0.5]{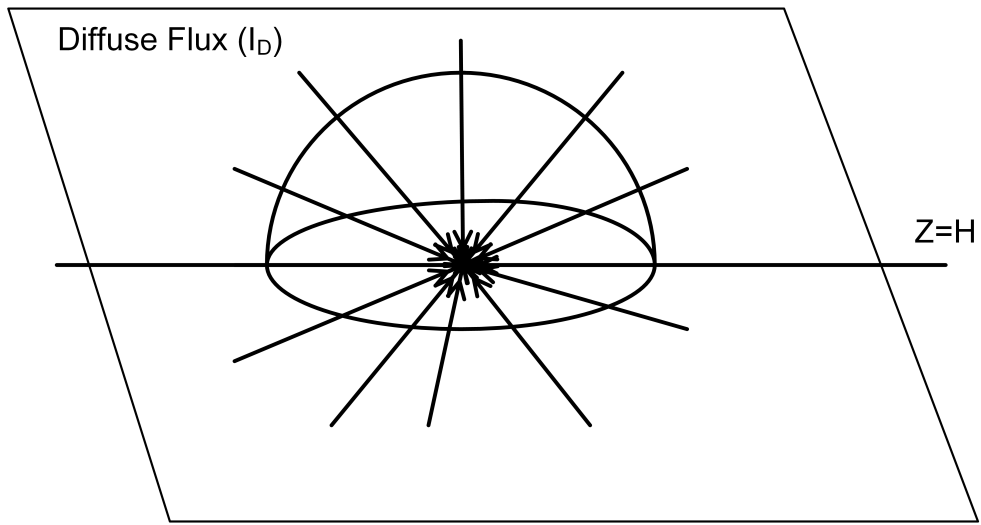}
\end{center}
\caption{A typical example of a downwelling irradiation}
\label{downwelling}
\end{figure}

The structure of the article is as follows. A model of phototaxis with absorption and scattering with the effects of diffuse and oblique collimated irradiation is formulated. The governing equations and boundary conditions are discussed next. This section is followed by deriving the solution for basic steady state, and then the linear stability equations are derived. Neutral stability curves for the problem are obtained numerically and  finally, the physical interpretation of the results are discussed. Afterwards, a comparison with the up-swimming phototaxis model has been made and the novelty of the proposed model is addressed thereafter.

\begin{figure}[!h]
\begin{center}
\includegraphics[scale=0.5]{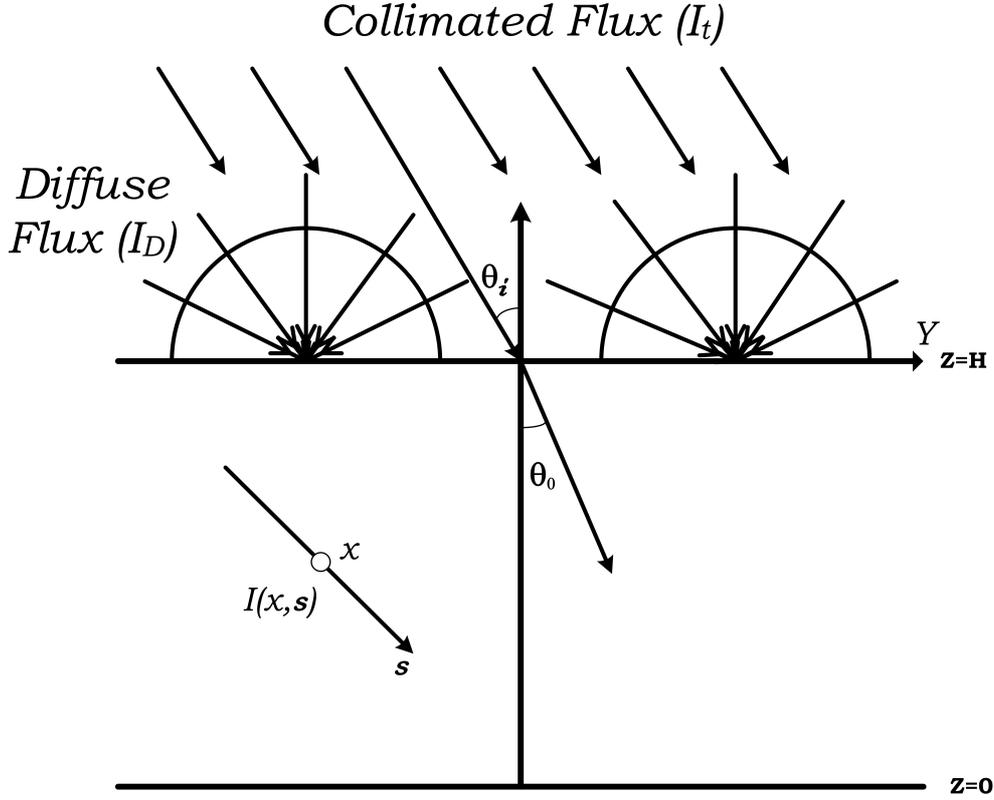}
\end{center}
\caption{Diffuse and oblique collimated irradiation incident on the surface of a suspension.}
\label{fig1}
\end{figure}

\section{GEOMETRY OF THE PROBLEM}
\label{secx2}

Consider the motion in a dilute suspension of phototactic algae within a layer of constant depth $H$ and infinite lateral extent [see Fig. \ref{fig3}]. It is assumed that the top and the bottom boundaries are to be non-reflecting. Here, the oblique collimated and diffuse irradiation both illuminate uniformly the algal suspension from above. The oblique collimated irradiation strikes the air-water interface  at a fixed off-normal angle $\theta_{i}$ and is transmitted across it at $z=H$ into water [see Fig. \ref{fig3}]. A rectangular cartesian coordinate system is choosen here and the $yz-$plane is the plane of incidence for the oblique collimated irradiation [see Fig. \ref{fig3}]. As the refractive index of the water is different from that of air, the angle of refraction $\theta_{0}$ is related to the angle of incidence $\theta_{i}$ by Snell's law, $\sin{\theta_{i}}=n_{0}\,\sin{\theta_{0}},$ where $n_{0}$ is the index of refraction of the water. The estimated value for refractive index of the water is $1.333$ approximately \cite{Daniel_1979} and the index of refraction of air has been assumed to be equal to unity.  The light incident on the algal suspension is absorbed and scattered thereafter, since the index of refraction of algae is not the same as that of water.

\section{PHOTOTAXIS WITH ABSORPTION AND SCATTERING}
\label{sec2}

Let $\mathrm{I}({\bmi x},{\bmi s})$ be the intensity of light at position ${\bmi x}$ and direction ${\bmi s}=\cos \theta \hat {\bmi z}+ \sin \theta (\cos \phi \hat {\bmi x} + \sin \phi \hat {\bmi y})$, where ${\bmi x}$ is measured relative to a rectangular cartesian coordinate system with the $z$-axis vertically up [see Fig. \ref{fig1}]. Here $\theta$ denotes the polar angle (measured from the $z$-axis) and $\phi$ denotes the azimuthal angle (measured in between the projection of the radiation intensity onto $xy$- plane and the $x$- axis) describing the unit vector ${\bmi s}$ (in spherical polar coordinate system).

The suspension is assumed to be illuminated uniformly by both the oblique collimated (unscattered) and diffuse (scattered) irradiation at the top (see Fig. \ref{fig1} and Sec. \ref{sec1} for justification). The radiative transfer equation (RTE) for an absorbing and scattering medium is given by \cite{modest}
\beq
\frac{\ud \mathrm{I}({\bmi x},{\bmi s})}{\ud s} + (a+\sigma_s)\mathrm{I}({\bmi x},{\bmi s}) = \frac{\sigma_s}{4\pi}\int_0^{4\pi}\mathrm{I}({\bmi x},{\bmi s}') p({\bmi s},{\bmi s}')\,d\Omega', \label{rads}
\eeq
where $s$ is the path lenth along ${\bmi s}$, $a$ is the absorption coefficient,  $\sigma_s$ is the scattering coefficient and $\Omega$ is the solid angle. The scattering phase function, or phase function, $p({\bmi s},{\bmi s}')$, is a probability density function that gives the angular distribution of light intensity scattered from direction ${\bmi s}'$ into certain other direction ${\bmi s}$. We have considered  here the case of isotropic scattering  (independent of direction) for simplicity and hence $p=1$.

Therefore, the light intensity on the top boundary surface at location ${\bmi x}_b$ is

$$
\mathrm{I}({\bmi x}_b=(x,y,H),{\bmi s}) = \mathrm{I}_t \, \delta({\bmi s}-{\bmi s}_0) +\dfrac{\mathrm{I_D}}{\pi},
$$

where $\delta$ is the Dirac-delta function, $\mathrm{I}_t$ is the magnitude of the oblique collimated irradiation in the direction ${\bmi s}_0=\cos{\left(\pi-\theta_{0}\right)}\hat{\bmi{z}}+\sin{\left(\pi-\theta_{0}\right)}\left(\cos{\phi_{0}}\hat{\bmi{x}}+\sin{\phi_{0}\hat{\bmi{y}}}\right)$ is the incident direction defined in spherical polar coordinates [Refer Panda \textit{et al.}\cite{PPS}]. Here, $\mathrm{I_D}$ is magnitude of the diffuse irradiation [refer Panda \textit{et al.} \cite{prmm16} for more details]. We shall get back to the phototaxis model of GPH \cite{gph} in the absence of diffuse  and oblique collimated irradiation (i.e. $\mathrm{I_D}=0$ and $\theta_{i}=0$) [see also Sec. \ref{steady} for justification on boundary intensity].

We assume that the absorption and scattering coefficient are propotional to the concentration. Thus $\sigma_s=\beta n({\bmi x})$ and $a = \alpha n({\bmi x})$ and the RTE for an isotropically scattering medium becomes

\beq
\frac{\ud \mathrm{I}({\bmi x},{\bmi s})}{\ud s} + (\alpha+\beta) n \mathrm{I}({\bmi x},{\bmi s}) = \frac{\beta n}{4\pi}\int_0^{4\pi}\mathrm{I}({\bmi x},{\bmi s}') \,d\Omega'. \label{rte1}
\eeq

The total intensity $G({\bmi x})$ at a point ${\bmi x}$ in the medium is 

$$
G({\bmi x})=\int_0^{4\pi}\mathrm{I}({\bmi x},{\bmi s})\,d\Omega,
$$ 
and the radiative heat flux ${\bmi q}({\bmi x})$ is 

\beq
{\bmi q}({\bmi x}) = \int_0^{4\pi}\mathrm{I}({\bmi x},{\bmi s})\,\left[ \left({\bmi s} \bcdot \hat{\bmi{x}} \right) \hat{\bmi{x}} + \left({\bmi s} \bcdot \hat{\bmi{y}} \right) \hat{\bmi{y}}+ \left({\bmi s} \bcdot \hat{\bmi{z}} \right) \hat{\bmi{z}} \right]\,d\Omega.
\label{flux}
\eeq 

Let ${\bmi p}$ be the unit vector corresponding to the swimming direction
and $<{\bmi p}>$ is the ensemble average of the swimming direction for all the 
cells in  an elemental volume. For many species of micro-organisms, the swimming 
speed is independent of illumination, position, time and direction \cite{hh:hh}
and we denote the ensemble-average swimming speed by $W_c$. The average 
swimming velocity is thus 

$$
{\bmi W}_c = W_c <{\bmi p}>.
$$ 

The mean swimming direction, $<{\bmi p}>$, is given by

\beq
<{\bmi p}> = -T(G) \frac{{\bmi q}}{\varpi+|{\bmi q}|}, \label{eq2}
\eeq

where $\varpi \ge 0$ is a constant and $T(G)$ is the phototaxis function 
such that

\beq
T(G)\quad
\left\{ \begin{array}{c}
\ge 0, \quad \mbox{if } G \le G_c, \\
< 0, \quad \mbox{if } G > G_c.\\
\end{array}
\right.
\eeq

$T(G)$ depends on the total light intensity $G$ reaching the cell. The exact functional form of $T$ will
depend on the species of micro-organisms \cite{vh:hv}. Here, the phototaxis function is composed using suitable trigonometric functions such that its plot will fit approximately to the corresponding phototactic response curve obtained from the experimental findings [see Fig. 1 of Vincent and Hill \cite{vh:hv}]. The negative sign in Eq. (\ref{eq2}) incorporates the fact that the source of light intensity for a micro-organism lies in the opposite direction to the radiative heat flux vector. Here, the mean swimming orientation is described by a multidirectional radiative heat flux rather than the unidirectional one as taken by Panda \textit{et al.} \cite{PPS}. The mean swimming direction is zero when the light intensity is critical ($T(G)=0$) or isotropic (${\bmi q}=0$). The constant $\varpi \ge 0$ is introduced to handle the case of isotropic light intensity. If the light intensity across the medium is not 
isotropic (such as the problem considered here), then we can take $\varpi=0$. Thus, the mean swimming direction in this article is given by Eq. (\ref{eq2}) with $\varpi=0$.

\section{THE CONTINUUM MODEL}
\label{sec3}

We assume a monodisperse cell population which can be modelled by a continuous distribution reminiscent to the previous models on bioconvection \cite{pk:kp}. The algal suspension is dilute so that the volume fraction of the cells is small and cell-cell interactions are negligible. Each cell has a volume $\vartheta$ and density $\rho+\Delta\rho$, where $\rho$ is the density of the fluid in which the cells swim and $\Delta\rho/\rho \ll 1$. ${\bmi u}$ is the average velocity of all the material in a small volume $\delta\!V$ and
$n$ is the cell concentration. Supposing that the suspension is incompressible, the average fluid velocity ${\bmi u}$ satisfies  

\begin{equation}
\bnabla\bcdot{\bmi u} = 0.  \label{eq3}
\end{equation}

We shall assume that the effect of cells on the suspension is dominated by
Stokeslets due to negative buoyancy and all other contributions of the cells 
to the bulk stress are sufficiently small to be neglected. Thus, neglecting 
all forces on the fluid except the cell's negative buoyancy, $n\,g\,\vartheta\,
\Delta\rho$ per unit volume where $g$ is the acceleration due to gravity, the
momentum equation under the Boussinesq approximation is

\begin{equation}
\frac{\mathrm{D}{\bmi u}}{\mathrm{D} t} = -\bnabla p_e + \mu{\nabla^2}{\bmi u} -
n\Delta\rho g \vartheta \hat {\bmi z}. \label{eq4}
\end{equation}

Here $\mathrm{D}/\mathrm{D}t=\p/\p t + {\bmi u}\bcdot\bnabla$ is the material
derivative, $p_e$ is the excess pressure above hydrostatic, $\hat {\bmi z}$ is 
a unit vector vertically upward, and $\mu$ is the viscosity of the suspension 
which is assumed to be that of the fluid. 

The equation for cell conservation is

\begin{equation}
\frac{\p n}{\p t} = -{\bmi\nabla}\bcdot{\bmi F}, \label{eq5}
\end{equation}
where  ${\bmi F}$ is the flux of cells. ${\bmi F}$ can be written as
\begin{equation}
{\bmi F} = n {\bmi u}+{n}W_c <{\bmi p}> - D {\bmi\nabla}{n}.  \label{eq6}
\end{equation}

Here the first term on the right-hand side of Eq. (\ref{eq6}) is the flux due to the advection of the cells by the bulk fluid flow and the second term arises due to the average swimming of the cells. The third term represents the random component of the cell locomotion. We choose the diffusivity tensor $\mathsf{\mathbf D}$ to be isotropic and constant and thus $\mathsf{\mathbf D}=D\mathsf{\mathbf I}$. Each cell is assumed to be purely phototactic and a spherical homogeneous body having uniform distribution of mass. Thus, the centre of mass of the body and its geometrical centre coincide. The cell flux vector in Eq. (\ref{eq6}) has been expressed via two important assumptions. First, phototaxis is in general a strong active orientation mechanism and the effect of viscous torque due to local fluid gradients, which might initiate horizontal shear, is neglected.  Second, the diffusion tensor is assumed to be constant, whereas it should be derived from the swimming velocity autocorrelation function. These assumptions allow us to remove the Fokker-Planck equation from the governing system for bioconvection. Thus, the resulting model is a valid limiting case in order to consider to understand the complexity of the problem before moving to more complex detailed model.

The lower boundary is taken as rigid while the upper boundary may be free or rigid
as in experiments. Indeed, even if the upper boundary is open to the air, cells
often collect at the top surface forming a rigid like packed layer. The boundary
conditions are

\begin{eqnarray}
{\bmi u}\bcdot {\hat {\bmi z}}&=&0\qquad \mbox{on}\quad z=0,H,\label{bcc1}\\
{\bmi F}\bcdot {\hat {\bmi z}}&=&0\qquad\mbox{on}\quad z=0,H.\label{bc2}
\end{eqnarray}

For rigid boundaries

\begin{equation}
{\bmi u}\times \hat {\bmi z}=0\qquad \mbox{on}\quad z=0,H, \label{bc3}
\end{equation}

while for a free boundary

\begin{equation}
\frac{\p^2}{\p z^2}({\bmi u}\bcdot\hat {\bmi z})=0. \label{bc4}
\end{equation}

We assume that the top boundary is exposed to uniform diffuse irradiation and oblique collimated radiation. We take the polar axis along the $z$ direction and thus, the boundary intensities are

%\begin{widetext}
\begin{subequations}
\label{bc56}
\begin{eqnarray}
\mathrm{I}(x,y,z=H,\theta,\phi) &=& \mathrm{I}_t \, \delta({\bmi s}-{\bmi s}_0)+\dfrac{\mathrm{I_D}}{\pi},\quad \pi/2\le\theta\le \pi, \label{bc5}\\
\mathrm{I}(x,y,z=0,\theta,\phi) &=& 0,\qquad  0\le\theta\le \pi/2.
\label{bc6}
\end{eqnarray}
\end{subequations}
%\end{widetext}

The governing equations are made dimensionless by scaling all lengths on $H$, the
depth of the layer, time on diffusive time scale $H^2/D$, and the bulk fluid
velocity on $D/H$. The appropriate scaling for the pressure is $\mu D/H^2$ and the
cell concentration is scaled on $\bar n$, the mean concentration. In terms of the
non-dimensional variables, the bioconvection equations become

\begin{equation}
\bnabla\bcdot{\bmi u} = 0,  \label{neq3}
\end{equation}

\begin{equation}
S_c^{-1}\frac{\mathrm{D}{\bmi u}}{\mathrm{D} t} = -\bnabla p_e + {\nabla^2}{\bmi u} -
R\,n \hat {\bmi z}, \label{neq4}
\end{equation}
\begin{equation}
\frac{\p n}{\p t} = -{\bmi\nabla}\bcdot{\bmi F}, \label{neq5}
\end{equation}

where

\begin{equation}
{\bmi F} = n {\bmi u}+{n}V_c <{\bmi p}> - {\bmi\nabla}{n}.  \label{neq6}
\end{equation}

Here $S_c=\mu/\rho D$ is the Schmidt number, $V_c$ is the scaled swimming speed, 
and $R$ is a Rayleigh number defined as

$$
R=\bar n \vartheta g \Delta\rho H^3/\nu D \rho.
$$

In dimensionless form, the boundary conditions  become

\begin{eqnarray}
{\bmi u}\bcdot {\hat {\bmi z}}&=&0\qquad \mbox{on}\quad z=0,1,\label{nbcc1}\\
{\bmi F}\bcdot {\hat {\bmi z}}&=&0\qquad\mbox{on}\quad z=0,1.\label{nbc2}
\end{eqnarray}

For rigid boundaries

\begin{equation}
{\bmi u}\times \hat {\bmi z}=0\qquad \mbox{on}\quad z=0,1, \label{nbc3}
\end{equation}

while for a free boundary

\begin{equation}
\frac{\p^2}{\p z^2}({\bmi u}\bcdot\hat {\bmi z})=0. \label{nbc4}
\end{equation}

The radiation transfer equation [see Eq. (\ref{rads})] becomes

\beq
\frac{\ud \mathrm{I}({\bmi x},{\bmi s})}{\ud s} + \kappa n \mathrm{I}({\bmi x},{\bmi s}) = \frac{\sigma n}{4\pi}\int_0^{4\pi}\mathrm{I}({\bmi x},{\bmi s}') \,d\Omega',  \label{rat2}
\eeq

where $\kappa=(\alpha+\beta) \bar n H$ is the nondimensional extinction coefficient and $\sigma=\beta \bar n H$ is the nondimensional scattering coefficient. The single scattering (scattering) albedo, a measure of the scattering efficiency of micro-organisms, is defined as $\omega=\dfrac{\sigma}{\kappa}.$ In terms  of scattering albedo $\omega$, Eq. (\ref{rat2}) can be written as:

\beq
\frac{\ud \mathrm{I}({\bmi x},{\bmi s})}{\ud s} + \kappa n \mathrm{I}({\bmi x},{\bmi s}) = \frac{\omega\kappa n}{4\pi}\int_0^{4\pi}\mathrm{I}({\bmi x},{\bmi s}') \,d\Omega'. \label{rat3}
\eeq

Here $\omega\in [0,1]$ and $\omega=0$ implies a purely absorbing (no scattering) medium whereas $\omega=1$ represents a purely scattering (no absorption) medium. In dimensionless form, the intensity at the top and bottom  becomes

%\begin{widetext}
\begin{subequations}
\label{nbc56}
\begin{eqnarray}
\mathrm{I}(x,y,z=1,\theta,\phi) &=& \mathrm{I}_t \, \delta({\bmi s}-{\bmi s}_0)+\dfrac{\mathrm{I_D}}{\pi},\quad \pi/2\le\theta\le \pi, \label{nbc5}\\
\mathrm{I}(x,y,z=0,\theta,\phi) &=& 0,\qquad  0\le\theta\le \pi/2.
\label{nbc6}
\end{eqnarray}
\end{subequations}
%\end{widetext}

\section{THE STEADY SOLUTION}
\label{steady}

\begin{figure}[!h]
\begin{center}
\includegraphics[scale=1]{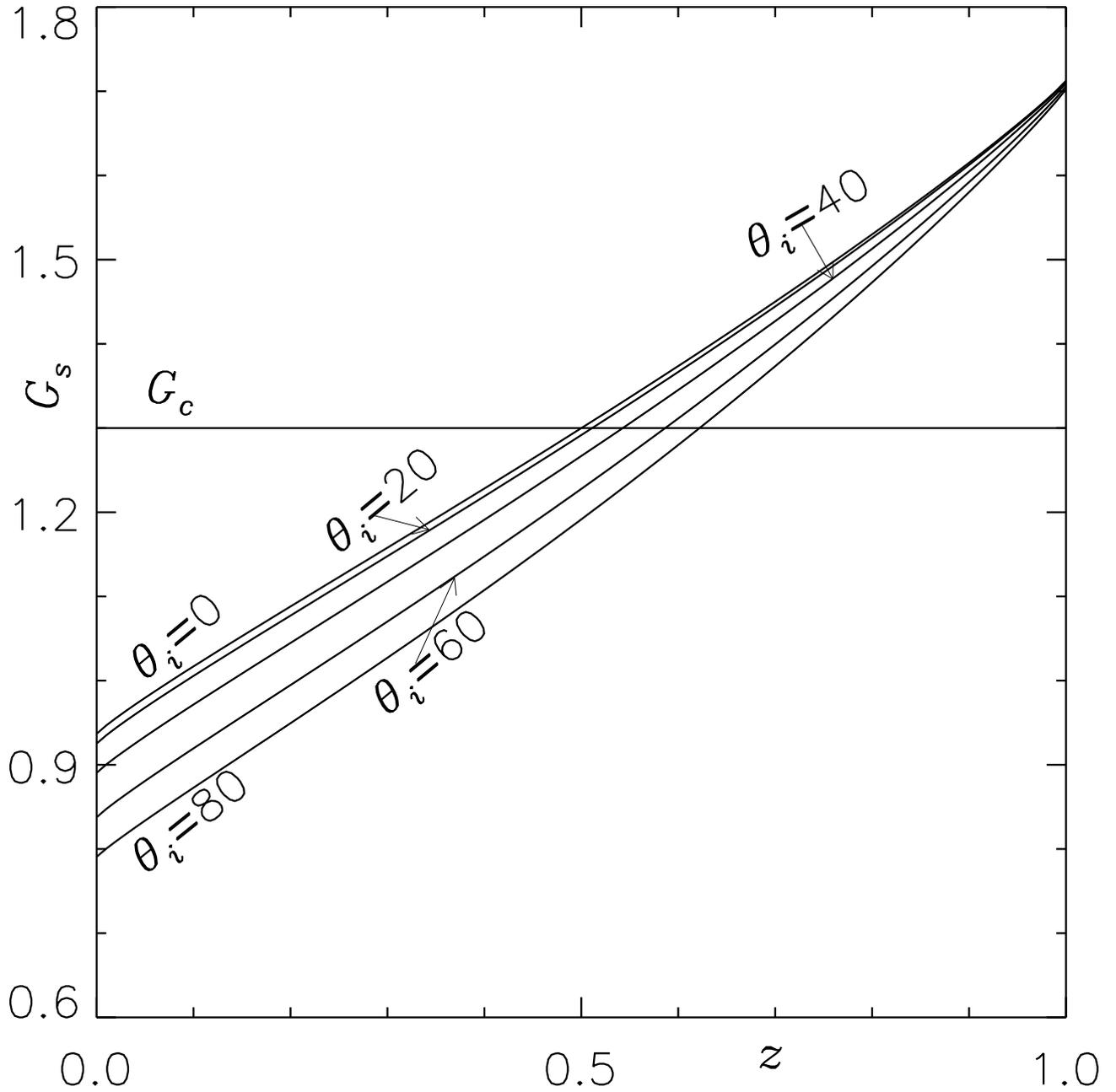}
\end{center}
\caption{Total intensity variation for different values of $\theta_{i}$ in a
uniform suspension. Here the other fixed parameters are $\kappa=0.5$, $\omega=0.4$, and $\mathrm{I_D}=0.26$.
}
\label{fig3}
\end{figure}

Equations (\ref{neq3})--(\ref{neq6}) and Eq. (\ref{rat3}) possess a static equilibrium 
solution in which

\beq
{\bmi u}=0,n=n_s(z)\quad\mbox{and}\quad \mathrm{I}=\mathrm{I}_s(z,\theta). \label{eqs}
\eeq

Then total intensity $G_s$ and radiative heat flux ${\bmi q}_s$ at the basic state are
given by

$$
G_s=\int_0^{4\pi}\mathrm{I}_s(z,\theta)\,d\Omega,\quad {\bmi q}_s = \int_0^{4\pi}\mathrm{I}_s(z,\theta)\,{\bmi s}\,d\Omega.
$$  
           
The equation governing $\mathrm{I}_s$ using direction cosines $(\xi=\sin{\theta} \cos{\phi},\eta=\sin{\theta} \sin{\phi},\nu=\cos{\theta})$ can 
be written as

\beq
\frac{d \mathrm{I}_s}{d z} + \frac{\kappa n_s \mathrm{I}_s}{\nu} =\frac{\omega\kappa
n_s}{4\pi\nu}G_s(z). \label{js} 
\eeq   
                   
We decompose the basic state intensity into a collimated part after attenuation, $\mathrm{I}_s^c$, and 
a diffused part which occurs due to scattering, $\mathrm{I}_s^d$, i.e. $\mathrm{I}_s = \mathrm{I}_s^c + \mathrm{I}_s^d$. The collimated part, $\mathrm{I}_s^c(z,\theta)$ is governed by

$$
\frac{d \mathrm{I}_s^c}{d z} + \frac{\kappa n_s \mathrm{I}_s^c}{\nu} =0,                 
$$  

subject to the boundary condition

$$
\mathrm{I}_s^c(1,\theta)=\mathrm{I}_t \, \delta({\bmi s}-{\bmi s}_0), \quad \pi/2\le\theta\le \pi.
$$

Now $\mathrm{I}_s^c$ is given by

$$
\mathrm{I}_s^c=\mathrm{I}_t \exp\left(\int_z^1\frac{\kappa
n_s(z')}{\nu}\,dz'\right)\delta({\bmi s}-{\bmi s}_0),
$$

and $\mathrm{I}_s^d$ is governed by

$$
\frac{d \mathrm{I}_s^d}{d z} + \frac{\kappa n_s \mathrm{I}_s^d}{\nu} =\frac{\omega\kappa
n_s}{4\pi\nu}G_s(z),                 
$$  
       
subject to the boundary conditions

\begin{subequations}
\label{snbc56}
\begin{eqnarray}
\mathrm{I}_s^d(1,\theta) &=& \dfrac{\mathrm{I_D}}{\pi},\, \quad \pi/2\le\theta\le \pi, \label{snbc5}\\
\mathrm{I}_s^d(0,\theta) &=& 0,\qquad  0\le\theta\le \pi/2.
\label{snbc6}
\end{eqnarray}
\end{subequations}

Eq.~(\ref{snbc5}) can be justified on the assumption that the incident radiation is diffuse on the azimuthal symmetric top surface ($z=1$) of the suspension and thus, $\mathrm{I}_s^d(1,\theta),\, \pi/2\le\theta\le \pi,$ is assumed to be direction independent and hence constant. Hence the magnitude of the diffuse irradiation, $\mathrm{I_D}$ reduces to $$\mathrm{I_D}=\pi \, \mathrm{I}_s^d(1,\theta).$$

Now the total intensity is decomposed as the sum of a collimated and diffused one i.e.
\beq
G_s=G_s^c+G_s^d, \label{Tintensity}
\eeq
 where
\begin{eqnarray*}
G_s^c&=&\int_0^{4\pi}\mathrm{I}_s^c(z,\theta)\,d\Omega=\mathrm{I}_t\exp\left(\dfrac{-\kappa\int_z^1
n_s(z')\,dz'}{\cos{\theta_{0}}}\right),\\
G_s^d &=&\int_0^{4\pi}\mathrm{I}_s^d(z,\theta)\,d\Omega. 
\end{eqnarray*}

We get the well known  Lambert-Beer law \cite{herdan} $G_s=G_s^c,$ for no
scattering. If we define the optical thickness as

$$
\tau=\int_z^1 \kappa n_s(z')\,dz',
$$ 

then the total intensity $G_s$ becomes a function of $\tau$ only. Further,
nondimensional total intensity, $\Lambda(\tau)=G_s(\tau)/\mathrm{I}_t$, satisfies the following Fredholm integral equation:\cite{modest,crosbie_82}
%\begin{widetext}
\beq
\Lambda(\tau) =  \frac{\omega}{2}\int_0^{\kappa}\Lambda(\tau')\,E_1\left(|\tau-\tau'|\right)\,d\tau' + e^{-\tau/\cos{\theta_{0}}} + 2\,\mathrm{I}_D \, E_2\left(\tau \right), \label{rts}
\eeq
%\end{widetext}
where $E_n(x)$ is the exponential integral of order $n$ \cite{chandras}. Eq. (\ref{rts}) is solved using the method of substraction of singularity \cite{press}.

The basic state radiative heat flux is written as

%\begin{widetext}
\begin{equation*}
{\bmi q}_s = \int_0^{4\pi}(\mathrm{I}_s^c+\mathrm{I}_s^d)\,{\bmi
s}\,d\Omega={\bmi q}_s^{c}+{\bmi q}_s^{d}=-\cos{\theta_{0}}\,\mathrm{I}_t\exp\left(\dfrac{-\kappa\int_z^1 n_s(z')\,dz'}{\cos{\theta_{0}}}\right)\hat{\bmi z} +
\int_0^{4\pi}\mathrm{I}_s^d(z,\theta)\,{\bmi s}\,d\Omega.
\end{equation*}    
%\end{widetext}

Since $\mathrm{I}_s^d(z,\theta)$ is independent of $\phi$, the $x$ and $y$ components
of ${\bmi q}_s$ vanish. Therefore, ${\bmi q}_s=-q_s\hat{\bmi z}$, where
$q_s=|{\bmi q}_s|$. The mean swimming direction becomes 
$$
<{\bmi p}_s>=-T_s\frac{{\bmi q}_s}{q_s}=T_s\hat{\bmi z},
$$

where $T_s=T(G_s)$.
 
The concentration $n_s(z)$ satisfies 

\begin{equation}
\frac{\ud n_s}{\ud z} = V_c\,T_s\,n_s, \label{eq13}
\end{equation}
which is supplemented by the cell conservation relation
\begin{equation}
\int_{0}^{1} n_s\,\ud z = 1.  \label{cellcon}
\end{equation}

Equations (\ref{rts})--(\ref{cellcon}) constitute a boundary value problem which 
is solved numerically using a shooting method. 

\begin{figure}[!h]
\begin{center}
\includegraphics[scale=1]{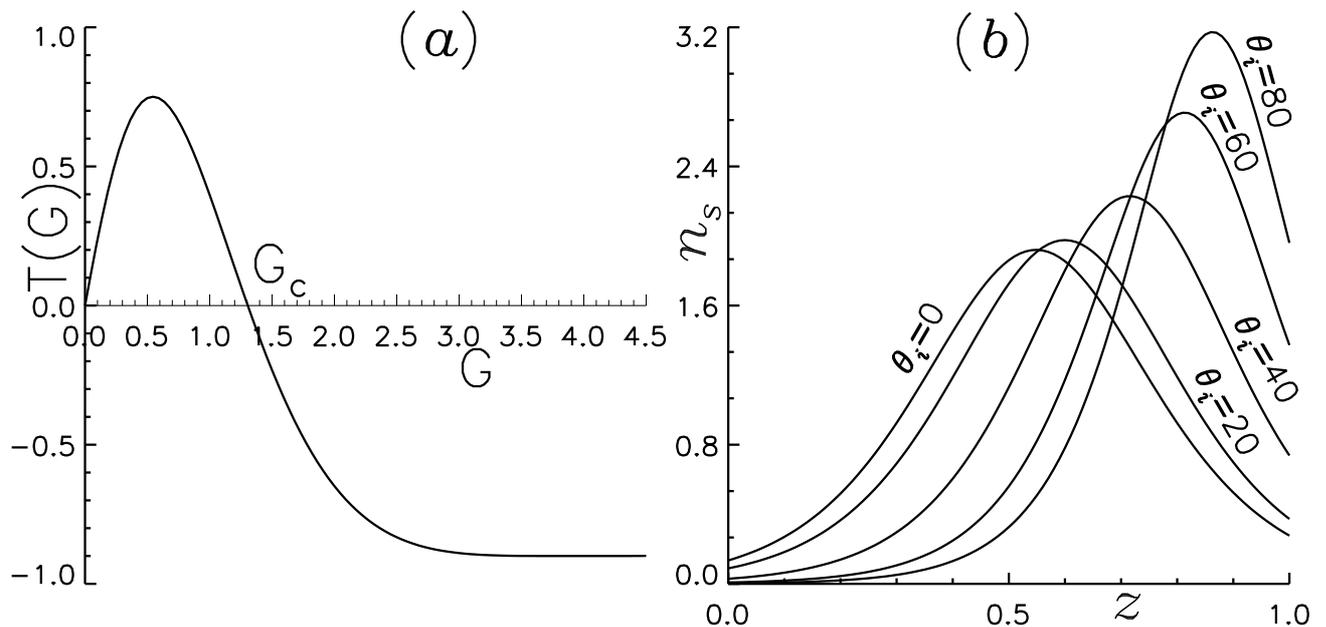}
\end{center}
\caption{(a) The phototaxis function with critical intensity $G_c=1.3$ [refer Eq.~(\ref{taxis1})] and (b) the corresponding base concentration profiles for different values of angle of incidence $\theta_{i}$. Fixed parameter values are $S_c = 20$, $V_c=15$, $\kappa=0.5$, $\omega=0.4$ and $\mathrm{I}_t=1$ and $\mathrm{I_D}=0.26$.}
\label{fig4}
\end{figure}

\begin{figure}[!h]
\begin{center}
\includegraphics[scale=1]{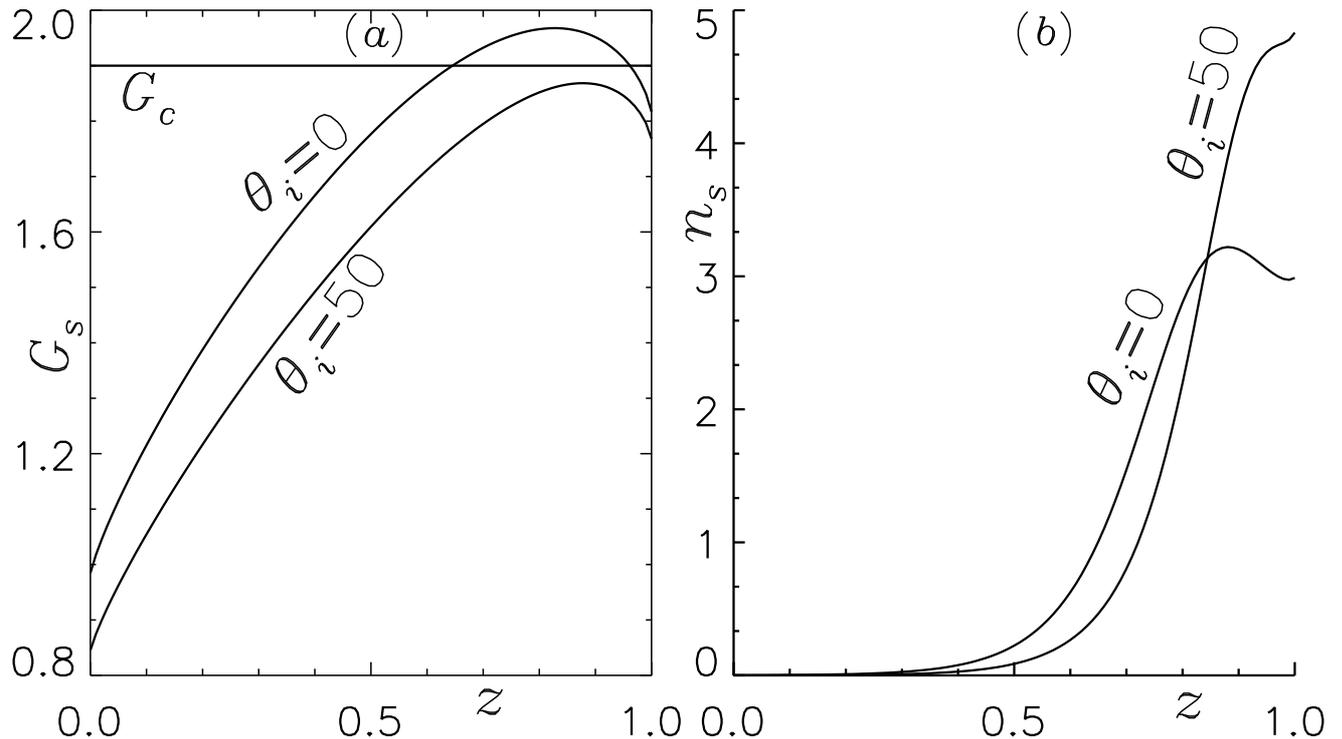}
\end{center}
\caption{(a) Variation of total intensity for $\theta_{i}=0,50$  and  in a uniform suspension (here the other fixed parameters are  $\kappa=1$ and $\omega=1$ and $\mathrm{I_D}=0.02$.) and (b) the corresponding base concentration profiles (here the other fixed parameters are $V_c=15$, $\kappa=1$, $\omega=1$ and $\mathrm{I}_t=1$). The phototaxis function used here has critical intensity $G_c=1.9$ [see Eq.~(\ref{taxis2})].}
\label{fig5}
\end{figure}

To estimate the parameters required for the present study, we assume that we are dealing with a phototactic micro-organism similar to \textit{Chlamydomonas}.  To retain the resulting model to be more rational and it can be comparable with earlier studies on (phototactic) bioconvection, we use the same parameter values as taken by \cite{vh:hv,HPK,gp,PPS} [see Table \ref{parameters}]. The range of the declination angle $\theta_{0}$ is restricted such that $0.661 \leq \cos{(\theta_{0})} \leq 1.0,$ and this implies that $0 \leq \theta_{0}(\textnormal{deg}) \leq 48.6$ approximately \cite{Daniel_1979}. We have estimated the approximate range of the angle of incidence $\theta_{i}$ as $0 \leq \theta_{i}(\textnormal{deg}) \leq 80$ for our proposed model [refer Daniel \textit{et al.}\cite{Daniel_1979}]. The radiation characteristics required for the present study are calculated as given in Ghorai and Panda~\cite{gp}. Thus, the optical depth $\kappa$ varies in the range from $0.25$ to $1$ for a $0.5$ cm deep suspension. The calulated scaled swimming speed for a $0.5$ cm and $1.0$ cm deep suspension are $V_c = 10$ and $V_c = 20$ respectively [see Table~\ref{parameters}]. The magnitude of the diffuse irradiation, $\mathrm{I_D}$ varies depending upon the prevailing overcast sky conditions and thus, $\mathrm{I_D}$ varies in the range from $0$ to $1.$

\begin{table}[!h]
\begin{center}
\caption{Estimates of typical parameters for a suspension of \textit{Chlamydomonas nivalis}.\label{parameters}}
\begin{tabular}{c c c}
\hline \hline
Cell radius         \quad  & $a$          \quad   & $10^{-3} \, \textnormal{cm}$ \\
Cell Volume         \quad  & $\vartheta$   \quad & $5 \times 10^{-10}\, \textnormal{cm}^3$ \\
Cell density ratio  \quad  & $\Delta{\rho}/\rho$ \quad & $5 \times 10^{-2}$ \\
Cell diffusivity    \quad  & $D$ \quad & $5 \times 10^{-4}\, \textnormal{cm}^2 \, \textnormal{s}^{-1}$\\
Swimming speed      \quad  & $W_c$    \quad & $10^{-2}\, \textnormal{cm} \, \textnormal{s}^{-1}$ \\
Mean concentration   \quad  & $\bar n$  \quad & $10^{6}\, \textnormal{cm}^{-3}$ \\
Kinematic viscosity  \quad   & $\mu/\rho$  \quad & $10^{-2}\,  \textnormal{cm}^2 \, \textnormal{s}^{-1}$\\
Schmidt number       \quad   & $S_c=\mu/\rho D$ \quad & $20$ \\
Scaled swimming speed \quad  & $V_c=W_c H/D$  \quad & $20H$\\
\hline \hline
\end{tabular}
\end{center} 
\end{table}

The magnitude of the oblique collimated irradiation at the top, $\mathrm{I}_t$ has been assumed to be equal to unity.  Figure~\ref{fig3} shows the variation of the total intensity, $G_s$, across the layer of a uniform suspension ($n=1$) for $\kappa=0.5$, $\omega=0.4$, $\mathrm{I_D}=0.26$ and different values of $\theta_{i}$(deg). 

A phototaxis function with critical intensity $G_c=1.3$ is considered here [see Fig.~\ref{fig4}(a)], whose mathematical form is given by

%\begin{widetext}
\begin{equation}
T(G)=0.8\, \sin{\left(\dfrac{3\,\pi}{2} \Xi(G)\right)}-0.1\, \sin{\left(\dfrac{\pi}{2}\, \Xi(G)\right)}, \quad \Xi(G)=\dfrac{1}{3.8}\, G\, \exp{\left[0.252\left(3.8-G\right)\right]}.
\label{taxis1}
\end{equation}
%\end{widetext}

For $0 \le$ $\theta_{i}$(deg) $\le 80,$ $G_s$ is monotonically decreasing across the suspension (see Fig.~\ref{fig4}). Now, consider the case when $V_c=15$, $\kappa=0.5$,  $\mathrm{I_D}=0.26$ and $\omega=0.4.$ In this case, when $\theta_{i}=0$, i.e., the critical intensity occurs around the mid-height of suspension, the cells accumulate around the mid-height of the domain. As $\theta_{i}$ increases to a higher non-zero value different from $0$, the maximum concentration increases and the location of the maximum concentration shifts towards the top of the domain. The maximum concentration is smallest when the maximum is located around the mid-height of the
domain [see Fig.~\ref{fig4}(b)]. The effects of the angle of incidence $\theta_{i}$ on basic concentration remain qualitatively similar for $0< \omega < 0.7$. Because the uniform total intensity decreases monotonically across the suspension for $0< \omega < 0.7$. 

Next, consider the case when $\omega = 1$. In this case $G_s$ does not decrease monotonically across the chamber as $\theta_{i}$ increases. Fig.~\ref{fig5} shows the taxis function

%\begin{widetext}
\begin{equation}
T(G)=0.8\, \sin{\left(\dfrac{3\,\pi}{2} \Xi(G)\right)}-0.1\, \sin{\left(\dfrac{\pi}{2}\, \Xi(G)\right)}, \quad \Xi(G)=\dfrac{1}{3.8}\, G\, \exp{\left[0.135\left(3.8-G\right)\right]}
\label{taxis2}
\end{equation}
%\end{widetext}

with critical intensity $G_c=1.9.$

Here $G_c=1.9$ occurs at two different depths, $z=0.95$  and $z=0.65$ in a uniform suspension, while keeping the parameters $\theta_{i}=0$, $\kappa=1$, $\omega=1$ and $\mathrm{I_D}=0.02$ as fixed  [see Fig.~\ref{fig5}(a)]. Thus, the cells above $z=0.95$ and below $z=0.65$ are positively phototactic, and those in between are negatively phototactic. As a result, the cells accumulate at around the top as well as near $z = 0.85$ in the basic steady state [see Fig.~\ref{fig5}(b)]. If we take $\theta_{i}=50,$  the positive phototaxis which occurs near the top of the suspension transits to negative phototaxis and the negative phototaxis in between $z=0.65$ and $z = 0.85$ transits to positive phototaxis. As a result, the cells accumulate around the top (i.e. single location) in the basic steady state [see Fig.~\ref{fig5}(b)].

\section{THE LINEAR STABILITY PROBLEM}

We consider a small perturbation of amplitude $\eps$ ($0<\eps\ll 1$) to the
equilibrium state (\ref{eqs}) so that

%\begin{widetext}
%\beq
$
{\bmi u}=\eps {\bmi u}_1+\mathcal{O}(\epsilon^2),\,n=n_s+\eps n_1+\mathcal{O}(\epsilon^2),\,\mathrm{I}=\left(\mathrm{I}_{s}^{c}+\mathrm{I}_{s}^{d}\right) + \eps \left(\mathrm{I}_{1}^{c}+\mathrm{I}_{1}^{d}\right)+\mathcal{O}(\epsilon^2) \quad\mbox{and}\quad <{\bmi p}> = <{\bmi p}_s> + \eps <{\bmi p}_1>+\mathcal{O}(\epsilon^2).
$
%\eeq
%\end{widetext}

If ${\bmi u}_1=(u_1,v_1,w_1),$ then, substituting the perturbed variables into
Eqs. (\ref{neq3})--(\ref{neq5}) and linearizing about the equilibrium state by collecting $\mathcal{O}(\eps)$ term, gives
\beq
\bnabla\bcdot{\bmi u}_1=0,\label{p1}
\eeq

\beq
S_c^{-1} \frac{\p {\bmi u}_1}{\p t} = -\bnabla p_e - R n_1 \hat{\bmi z} + \nabla^2{\bmi u}_1, \label{p2}
\eeq

\beq
\frac{\p n_1}{\p t} + w_1 \frac{\ud n_s}{\ud z} +V_c\,\bnabla \bcdot(<{\bmi
p}_s> n_1 + <{\bmi p}_1> n_s)=\nabla^2 n_1. \label{p3}
\eeq

If $G=G_s+\eps G_1+\mathcal{O}(\epsilon^2)=\left(G_s^c+ \eps G_{1}^{c}\right)+ \left(G_s^d+\eps G_{1}^{d}\right)+\mathcal{O}(\epsilon^2)$, then the steady collimated total intensity is perturbed as $\mathrm{I}_t\exp\left(\dfrac{-\kappa\int_z^1 \left(n_s(z')+\eps n_1+\mathcal{O}(\epsilon^2)\right)\,dz'}{\cos{\theta_{0}}}\right)$ and after simplification, we get 
\beq
G_{1}^{c}=\mathrm{I}_t\exp\left(\dfrac{-\kappa \int_z^1 n_s(z')\,dz'}{\cos{\theta_{0}}}\right) \left(\dfrac{\kappa \int_{1}^{z} n_{1} dz'}{\cos{\theta_{0}}}\right)\label{cr1}
\eeq

Similarly, $G_{1}^{d}$ is given by
\beq
G_{1}^{d} = \int_0^{4\pi}\mathrm{I}_{1}^{d}({\bmi x},{\bmi s})\,d\Omega.\label{r1}
\eeq
Similarly, for the radiative heat flux ${\bmi q}={\bmi q}_s+\eps {\bmi q}_1+\mathcal{O}(\epsilon^2)=\left({\bmi q}_{s}^{c}+{\bmi q}_{s}^{d}\right)+\eps \left({\bmi q}_{1}^{c}+{\bmi q}_{1}^{d}\right)+\mathcal{O}(\epsilon^2)$, we  find

\begin{eqnarray*}
{\bmi q}_{1}^{c} &=& -\mathrm{I}_t\exp\left(\dfrac{-\kappa \int_z^1 n_s(z')\,dz'}{\cos{\theta_{0}}}\right) \left(\dfrac{\kappa \int_{1}^{z} n_{1} dz'}{\cos{\theta_{0}}}\right) \cos{\theta_{0}} \hat{\bmi z}, \\
{\bmi q}_{1}^{d} &=& \int_0^{4\pi}\mathrm{I}_{1}^{d}({\bmi x},{\bmi s})\,{\bmi s}\,d\Omega.\label{r2}
\end{eqnarray*}

Now the expression
$$
-T(G_s+\eps G_1+\mathcal{O}(\epsilon^2))\frac{{\bmi q}_s + \eps {\bmi q}_1+\mathcal{O}(\epsilon^2)}{|{\bmi q}_s + \eps {\bmi
q}_1+\mathcal{O}(\epsilon^2)|}-T_s\hat{\bmi z},
$$
on simplification at $\mathcal{O}(\eps)$, gives perturbed swimming direction
\beq
<{\bmi p}_1>=G_1\frac{\ud T_s}{\ud G} \hat{\bmi z}-T_s\frac{{\bmi q}_1^H}{q_s}. \label{pp}
\eeq
Please note that the second term of the right hand side of the equality sign of the above expression i.e. $T_s\frac{{\bmi q}_1^H}{q_s}$ becomes zero for a non-scattering algal suspension illuminated by oblique collimated irradiation [refer Panda \textit{et al.} \cite{PPS} for details] and here ${\bmi q}_1^H$ is the horizontal component of perturbed radiative heat flux ${\bmi q}_1$. 

Substituting Eq. (\ref{pp}) into Eq. (\ref{p3}) and simplifying we get
%\begin{widetext}
\beq
\frac{\p n_1}{\p t} + w_1 \frac{\ud n_s}{\ud z} + V_c\frac{\p}{\p
z}\left(T_s n_1 + n_s \frac{\ud T_s}{\ud G} G_1\right) - V_c n_s
\frac{T_s}{q_s}\left(\frac{\p q_1^x}{\p x} + \frac{\p q_1^y}{\p y}\right)=\nabla^2 n_1. \label{p4}
\eeq
%\end{widetext}
By elimination of $p_e$ and horizontal component of ${\bmi u}_1$, Eqs. (\ref{p1}), (\ref{p2}) 
and Eq. (\ref{p4}) can be reduced to two equations for $w_1$ and $n_1$. These quantities can 
then be decomposed into normal modes such that
\beq
w_1=W(z) \exp(\gamma t + i(lx+my)),\quad n_1=\Theta(z)\exp(\gamma t + i(lx+my)). \label{q1}
\eeq
The equation governing diffuse perturbed intensity $\mathrm{I}_{1}^{d}$ using direction cosines $(\xi,\eta,\nu)$ can be written as
%\begin{widetext}
\beq
\xi \frac{\p \mathrm{I}_{1}^{d}}{\p x} + \eta \frac{\p \mathrm{I}_{1}^{d}}{\p y} + \nu \frac{\p \mathrm{I}_{1}^{d}}{\p
z} + \kappa n_s \mathrm{I}_{1}^{d} =\frac{\omega\kappa}{4\pi}(n_s G_{1}^{c}+ n_s G_{1}^{d}+ G_s n_1)-\kappa \mathrm{I}_s n_1, \label{q2} 
\eeq  
%\end{widetext}
subject to the boundary condition
\begin{subequations}
\label{j1}
\begin{eqnarray}
\mathrm{I}_{1}^{d}(x,y,1,\xi,\eta,\nu)=0,\,\pi/2\le \theta\le \pi,0\le\phi\le 2\pi,
\label{j1a}\\
\mathrm{I}_{1}^{d}(x,y,0,\xi,\eta,\nu)=0,\, 0\le \theta\le \pi/2,0\le\phi\le 2\pi.
\label{j1b}
\end{eqnarray} 
\end{subequations}
The form of Eq. (\ref{q2}) suggests the following expression for $\mathrm{I}_{1}^{d}$:
$$
\mathrm{I}_{1}^{d}=\Psi^{d}(z,\xi,\eta,\nu) \exp(\gamma t + i(lx+my)).
$$
From Eqs. (\ref{cr1}) and (\ref{r1}) we get
%\begin{widetext}
\begin{eqnarray}
G_{1}^{c}& =&\left[\mathrm{I}_t\exp\left(\dfrac{-\kappa \int_z^1 n_s(z')\,dz'}{\cos{\theta_{0}}}\right) \left(\dfrac{\kappa \int_{1}^{z} \Theta(z') dz'}{\cos{\theta_{0}}}\right)\right]\exp(\gamma t + i(lx+my))\nonumber\\
&=&\mathcal{G}^{c}(z) \,\exp(\gamma t + i(lx+my)),\label{pcollimated}\\
G_{1}^{d} &=& \mathcal{G}^{d}(z)\exp(\gamma t + i(lx+my))\nonumber\\
&=&\left(\int_0^{4\pi}\Psi^{d}(z,\xi,\eta,\nu)\,d\Omega\right)\exp(\gamma t + i(lx+my)) \label{pdiffused},
\end{eqnarray}
%\end{widetext}
where
$\mathcal{G}(z)=\mathcal{G}^{c}(z)+\mathcal{G}^{d}(z)$. Note that $\mathcal{G}^c$ given in Eq. (\ref{pcollimated}) is the perturbed intensity given in Ref.\cite{PPS} for no scattering. 

Now $\Psi^{d}$ satisfies
%\begin{widetext}
\beq
\frac{\ud \Psi^{d}}{\ud z} + \frac{[i(l\xi+m\eta)+\kappa n_s]}{\nu}\Psi^{d} =
\frac{\omega \kappa}{4\pi\nu}[n_s \mathcal{G}^{c}+n_s \mathcal{G}^{d} + G_s\Theta]-\frac{\kappa}{\nu}\mathrm{I}_s\Theta, \label{q5}
\eeq
%\end{widetext}
subject to the boundary condition
\begin{subequations}
\label{q5e}
\begin{eqnarray} 
\Psi^{d}(1,\xi,\eta,\nu)=0,\, \pi/2\le \theta\le \pi,0\le\phi\le
2\pi,\label{q5a}\\
\Psi^{d}(0,\xi,\eta,\nu)=0,\, 0\le \theta\le \pi/2,0\le\phi\le
2\pi.\label{q5b}    
\end{eqnarray} 
\end{subequations}

Equation (\ref{q5}) is an integro-differential equation which is solved using iteration.

Similarly from Eq. (\ref{r2}), we have
$$
{\bmi q}_1^H=[q_1^x,q_1^y]=\left[P(z),Q(z)\right] \exp(\gamma t + i(lx+my)), 
$$
where 
$$
P(z) = \int_0^{4\pi}\Psi^{d}(z,\xi,\eta,\nu)\,\xi\,d\Omega,\quad Q(z) =
\int_0^{4\pi}\Psi^{d}(z,\xi,\eta,\nu)\,\eta\,d\Omega.
$$
The linear stability equations become
%\begin{widetext}
\beq
\left(\gamma S_c^{-1}+k^2-\frac{\ud^2}{\ud z^2}\right)\left(\frac{\ud^2}{\ud
z^2}-k^2\right)W=R k^2 \Theta,\label{lin1}
\eeq 
\beq
\left(\gamma+k^2-\frac{\ud^2}{\ud z^2}\right)\Theta + V_c \frac{\ud}{\ud
z}\left(T_s\Theta+n_s\frac{\ud T_s}{\ud G}\mathcal{G}\right)-i\frac{V_c n_s
T_s}{q_s}(lP+mQ) = -\frac{\ud n_s}{\ud z}W,\label{lin2} 
\eeq
%\end{widetext}
subject to the boundary conditions
\beq
\frac{\ud\Theta}{\ud z}-V_c T_s \Theta - V_c n_s \frac{\ud T_s}{\ud G} \mathcal{G}=0\qquad \mbox{at}\quad z=0,1, \label{linbc1}
\eeq
and for rigid boundaries
\beq
W=\frac{\ud W}{\ud z}=0 \qquad \mbox{at}\quad z=0,1. \label{linbc2}
\eeq
At a free surface the last condition in (\ref{linbc2}) is replaced by
\beq
\frac{\ud^2 W}{\ud z^2}=0. \label{linbc3} 
\eeq
Here $k=\sqrt{l^2+m^2}$, is the overall nondimensional wavenumber. Equations (\ref{lin1})--(\ref{linbc2}) form an eigen value problem for $\gamma$ as a functions of the dimensionless parameters  $\theta_{0}$(deg), $l,m,V_c,\kappa,\sigma$ and $R$. The basic state becomes unstable whenever $\mathrm{Re}(\gamma)>0$. Using notation $D={d}/{dz},$ Eq. (\ref{lin2}) becomes   
%\begin{widetext}
\beq
\Gamma_0(z)+\Gamma_1(z)\int_1^z \Theta\,dz +
\left[\gamma+k^2+\Gamma_2(z)\right]\Theta +
V_c T_s\,D\Theta - D^2\Theta =-(Dn_s)W, \label{linc}  
\eeq
%\end{widetext}
where
\begin{subequations}
\label{gam}
\begin{eqnarray}
\Gamma_0(z)&=&V_c\,D\left(n_s\frac{dT_s}{dG}\mathcal{G}^d\right)-i\frac{V_cn_sT_s}{q_s}(lP+mQ),\label{ga}\\
\Gamma_1(z)&=& \left(\dfrac{\kappa}{\cos{\theta_{0}}}\right) V_c D\left(n_s {G}_s^c \frac{dT_s}{dG}\right), \label{gb} \\
\Gamma_2(z)&=& 2\left(\dfrac{\kappa}{\cos{\theta_{0}}}\right) V_c n_s {G}_s^c \frac{dT_s}{dG}+V_c
\frac{dT_s}{dG}D{G}_s^d.\label{gc}
\end{eqnarray}
\end{subequations}
In the absence of scattering (i.e. $\omega=0$), $\Gamma_0(z)\equiv 0$,
${G}_s={G}_s^c,$ and the last term of Eq. (\ref{gc}) becomes zero. Then 
Eq. (\ref{linc}) becomes the same as that given in Ref. Panda \textit{et al.} \cite{PPS}. 

Introducing a new variable 
\beq
\Phi(z)=\int_1^z \Theta(z')\,dz', \label{extrav}
\eeq
the linear stability equations become
%\begin{widetext}
\beq
\left(\gamma S_c^{-1}+k^2-\frac{\ud^2}{\ud z^2}\right)\left(\frac{\ud^2}{\ud
z^2}-k^2\right)W=R k^2 D\Phi,\label{linwf}
\eeq     
\beq
\Gamma_0(z)+\Gamma_1(z)\Phi +
\left[\gamma+k^2+\Gamma_2(z)\right]D\Phi +
V_c T_s\,D^2\Phi - D^3\Phi =-(Dn_s)W. \label{lincf}  
\eeq  
%\end{widetext} 
The boundary conditions given in Eqs. (\ref{linbc2}) and (\ref{linbc3}) remain
the same except Eq. (\ref{linbc1}) which becomes
\beq
D^2\Phi-V_c T_s D\Phi - V_c n_s \frac{\ud T_s}{\ud G} \mathcal{G}=0\qquad \mbox{at}\quad
z=0,1, \label{linbc1f}
\eeq
and an extra boundary condition
\beq
\Phi(z)=0,\qquad \mbox{at}\quad z=0, \label{linbc2f}
\eeq
which is follows from Eq. (\ref{extrav}).

\section{SOLUTION PROCEDURE}

Numerical solutions to Eqs. (\ref{linwf}) and(\ref{lincf}) with appropriate boundary conditions are calculated with a fourth-order accurate, finite-difference scheme  based on Newton-Raphson-Kantorovich (NRK) iterations \cite{cm}. This scheme is used to calculate the neutral stability curves in the $(k,R)$-plane or the growth rate, Re$(\gamma)$, as a function of $R$ for a fixed set of other parameters. Initially, values of $S_c$, $V_c$, $\kappa$, $\omega$,  $k$, $\theta_{0}$ and $\mathrm{I_D}$ are supplied, and the values of $W$ and $\Phi$ are estimated either from the previous numerical results, or by imposing sinusoidal variation in $W$ and $\Phi$. Once a solution is obtained, this solution can be used as an initial guess for the neighboring parameter values.

For a given set of other parameter ranges, there are infinite number of branches of the neutral curve $R^{(n)}(k), \left(n = 1 , 2 , 3 , \cdots \right),$ each one representing a particular solution of the linear stability problem. The most interesting solution branch is the one on which $R$ has its minimum value $R_c$. The most unstable solution is recognized as the pair $\left(k_c, R_c\right),$ from which the wavelength of the initial disturbance may be obtained as $\lambda_{c}=2\pi/k_{c}$.  The bioconvective solutions consist of convection cells stacked one above another along the depth of the suspension. A solution is said to be of mode $n$ if it has n convection cells stacked vertically one on another. In many instances, the most unstable solution occurs on the $R^{(1)}(k)$ branch of the neutral curve and it is mode $1$.

A neutral curve is defined as the locus of points where Re$(\gamma) = 0$. If in addition Im$(\gamma) = 0$ on such a curve, then the principle of exchange of stabilities is said to be valid and the bioconvective solution is called stationary (non-oscillatory). Alternatively, if Im$(\gamma) \ne 0$ then oscillatory solutions exist. If the most unstable solution remains on the oscillatory branch of the neutral curve, then the solution is called overstable. When there is a competition between the stabilizing and destabilizing processes, oscillatory solution arises usually. When oscillatory solution occurs, a single oscillatory branch of the neutral curve meets the stationary branch of the neutral curve at $k = k_b$ and exists for $k \le k_b$. 

\begin{figure}[!h]
\begin{center}
\includegraphics[scale=1]{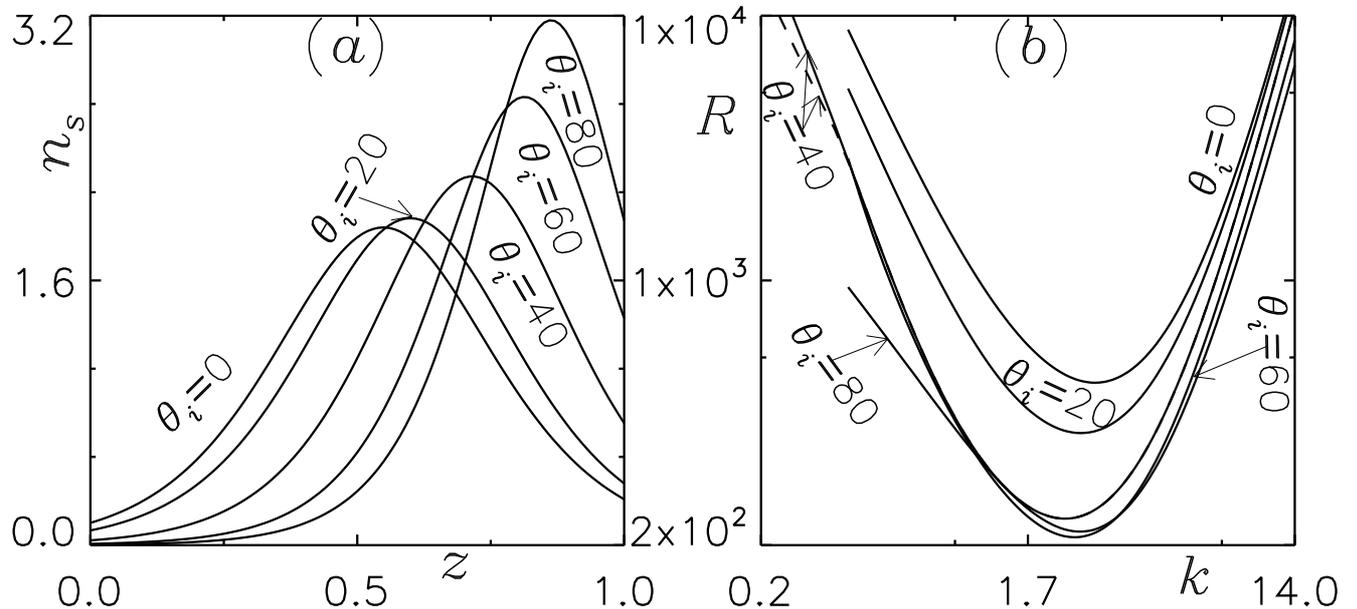}
\end{center}
\caption{(a) Base concentration profiles and (b) neutral curves as $\theta_{i}$ is increased. Fixed parameter values are $S_c = 20$, $V_c=15$, $\kappa=0.5$, $\omega=0.4$,$\mathrm{I_D}=0.26$ and $\mathrm{I}_t=1$.}
\label{fig6}
\end{figure}

\section{NUMERICAL RESULTS}

We have systematically investigated the effect of angle of incidence $\theta_{i}$(deg) by varying it between $0$ to $80$ (i.e. $0 \le \theta_{i} \le 80$),  keeping the other parameters $S_c$, $\mathrm{I}_t$, $G_c$, $V_c$, $\omega$, $\mathrm{I_D}$ and $\kappa$ fixed. Due to large number of parameter values, it is difficult to obtain a comprehensive picture across the whole parameter domain. Thus we take a discrete set of fixed parameter ranges which are of physical interest to study their effect on the onset of bioconvection. The values of $S_c = 20$, and $\mathrm{I}_t = 1$ are kept fixed throughout. The representative values of the cell swimming speed, the extinction coefficient and scattering albedo are $V_c =10, 15, 20$, $\kappa = 0.5, 1.0$, and $\omega \in [0:1]$  respectively. The value of the diffuse irradiation $\mathrm{I_D}$ is selected so that the maximum cell concentration in the basic state occurs around the mid-height ($z=1/2$) of the suspension. Based on self-shading by the algae, the results obtained for a discrete set of parameter ranges are divided into two categories. We address them separately via two cases.

\begin{figure}[!h]
\begin{center}
\includegraphics[scale=1]{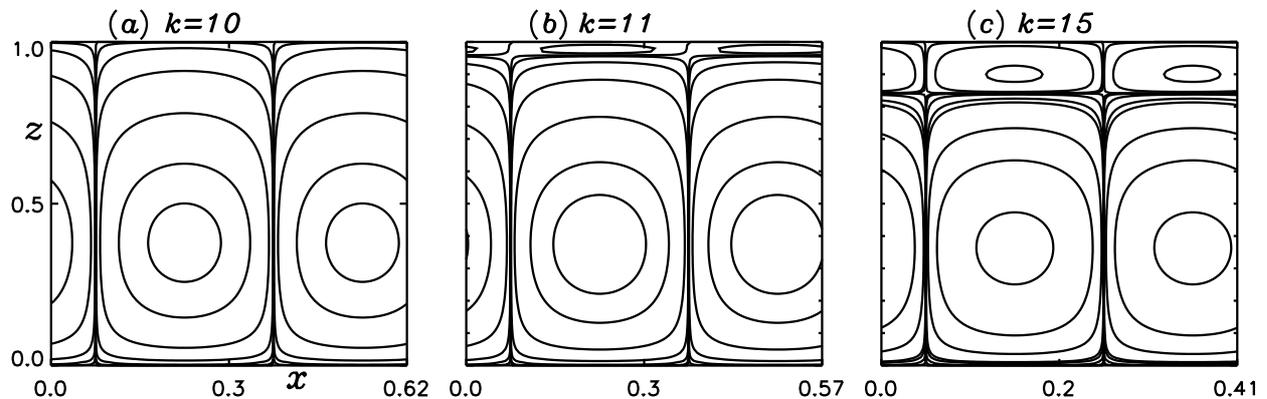}
\end{center}
\caption{Existence of mode $2$ solutions on $R^{(1)}(k)$ branch and the corresponding flow pattern of the perturbed velocity component (eigenmode) $w_1$. The small convection cells form near the top of the chamber. Fixed parameter values are $\mathrm{I_D}=0.26$, $V_c=15$, $\kappa=0.5$, $\omega=0.4$, $\theta_{i}=0$ and $\mathrm{I}_t=1$.}
\label{m2l}
\end{figure}

\subsection{Case-I : Weak-scattering algal suspension} 
\label{di}

\begin{figure}[!h]
\begin{center}
\includegraphics[scale=1]{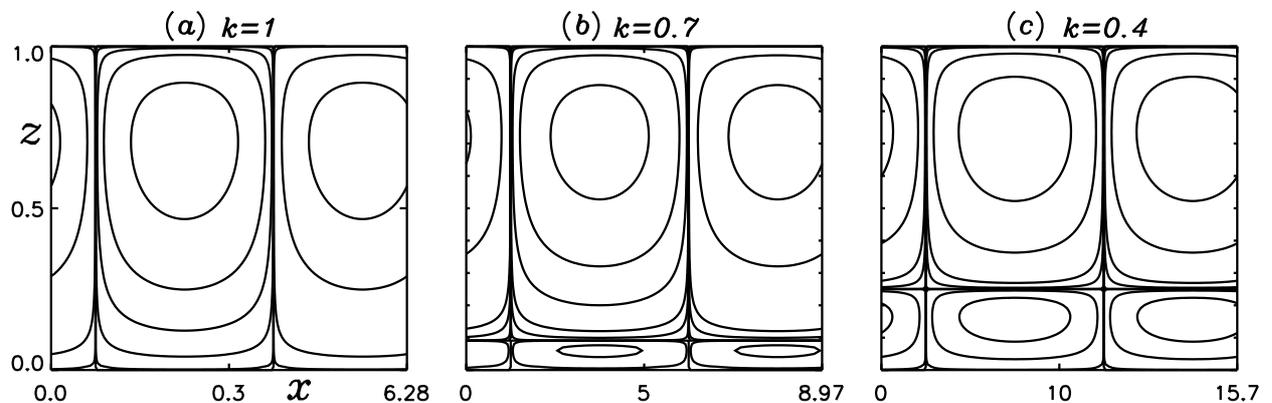}
\end{center}
\caption{Flow pattern of the perturbed velocity component (eigenmode) $w_1$. The small convection cells form near the bottom of the chamber. Fixed parameter values are $\mathrm{I_D}=0.26$, $V_c=15$, $\kappa=0.5$, $\omega=0.4$, $\theta_{i}=40$ and $\mathrm{I}_t=1$.}
\label{m2r}
\end{figure}

To study the effects of angle of incidence on bioconvection via a weak-scattering algal suspension, here we consider the case when the self-shading (absorption) is significant by selecting a lower value of scattering albedo $\omega$. Furthermore, self-shading (absorption) is more (less) effective when $\kappa=1$ ($\kappa=0.5$). Thus, the representative vaules of the governing parameter ranges taken are $V_c=10, 15, 20,$   $\kappa=0.5,1,$ $\mathrm{I_D}=0.26,$ and $\omega=0.4$ to include the effects of weak-scattering too. We take a phototaxis function here which has critical intensity $G_c=1.3$ [see Fig.~\ref{fig4}(a) and Eq.~(\ref{taxis1})].

\subsubsection{$V_c=15$} 
\label{}

\subsection*{Extinction coefficient $\kappa=0.5$} 
\label{ss}

\begin{figure}[!h]
\begin{center}
\includegraphics[scale=1]{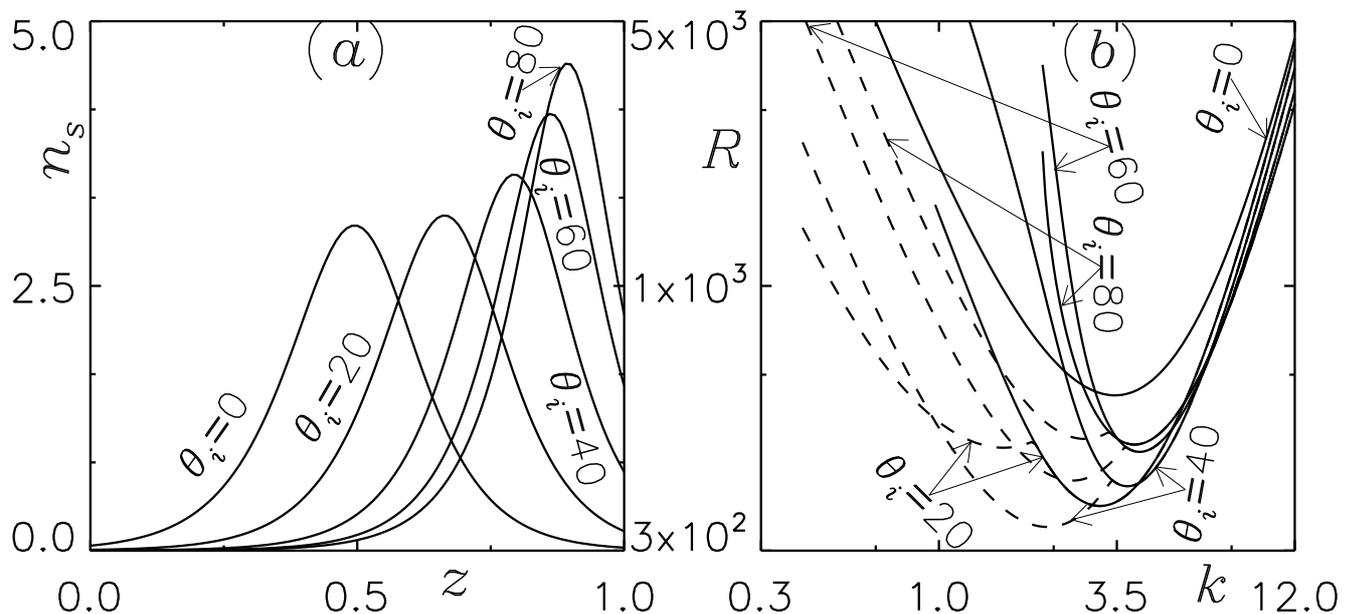}
\end{center}
\caption{(a) Base concentration profiles and (b) neutral curves showing stationary (solid lines) and oscillatory branches (dashed lines) as $\theta_{i}$ is increased. Fixed parameter values are $S_c = 20$, $V_c=15$, $\kappa=1$, $\omega=0.4$, $\mathrm{I_{D}}=0.5$ and 
$\mathrm{I}_{t}=1$.}
\label{fig7}
\end{figure}

The base concentration profiles and the corresponding neutral curves are shown in Fig.~\ref{fig6}, when the parameters $V_c=15$, $\kappa=0.5,$ $\mathrm{I_D}=0.26,$ and  $\omega=0.4$ are kept fixed and they serve to illustrate the effect of the angle of incidence for  $0 \le \theta_{i} \le 80$. For $\theta_{i}=0$, the base concentration develops maximum at around the mid-height of the suspension.
As $\theta_{i}$ is increased to $20$, the location of the maximum base concentration shifts towards the top of the suspension. Also, the width (thickness) of the upper stable layer overlying the unstable layer decreases as compared to the case when  $\theta_{i}=0$. As a result, the buoyancy of which tends to inhibit convective fluid motions becomes less effective and $R_c$ decreases. When $\theta_{i}$ is increased further to $40,$ the maximum concentration is located around three-quarter height of the suspension. In this case, a single oscillatory branch bifurcates from the stationary branch at wavenumber $k=0.42$ and defines a locus of points for $k \le 0.42.$ But, the most unstable solution still remains on the stationary branch. Also, both the critical Rayleigh number and wave number decrease. As $\theta_{i}$ becomes $60,$ the base concentration profile becomes steeper. This results in increment in critical wavenumber in comparison to the previous case, but the critical Rayleigh number decreases further. When $\theta_{i}=80,$ the steepness in the base concentration increases and width (thickness) of the upper stable layer decreases. As a result, both the critical Rayleigh number and wave number increase [see Fig.~\ref{fig6}].

It is also noticed that, along certain parts of the neutral curves, the bioconvective solution corresponding to $R^{(1)}(k)$ branch is mode $2$ not mode $1$. As the value of $k$ is increased above $k_c$, along such a branch, the single convection cell which extends throughout the depth of the layer becomes augmented by a second small convection cell, which originates at the top of the layer and grows in height as $k$ is increased further. This feature is observed on the $R^{(1)}(k)$ branch for $\theta_{i}=0$. But, the whole of the branch is mode $2$ here, except for a region $k < k_c$, where the solution is mode $1$. Thus, the solution of the linear stability problem is referred to as mode $2$ here [see Fig.~\ref{m2l}]. When $\theta_{i}=20,$ the bioconvective solution corresponding to $R^{(1)}(k)$ branch is mode $2$ too. As $\theta_{i}$ is increased further to $40,$ another feature is observed along certain parts of the neutral curves and is illustrated here. In such a case, if $k$ is decreased to a lower value than $k_c$ along the $R^{(1)}(k)$ branch, a second small convection cell appears at the bottom of the layer and grows in height as $k$ is decreased further.  Thus, the most unstable solution in this case is called mode $1$ [see Fig.~\ref{m2r}].

\begin{figure}[!h]
\begin{center}
\includegraphics[scale=1]{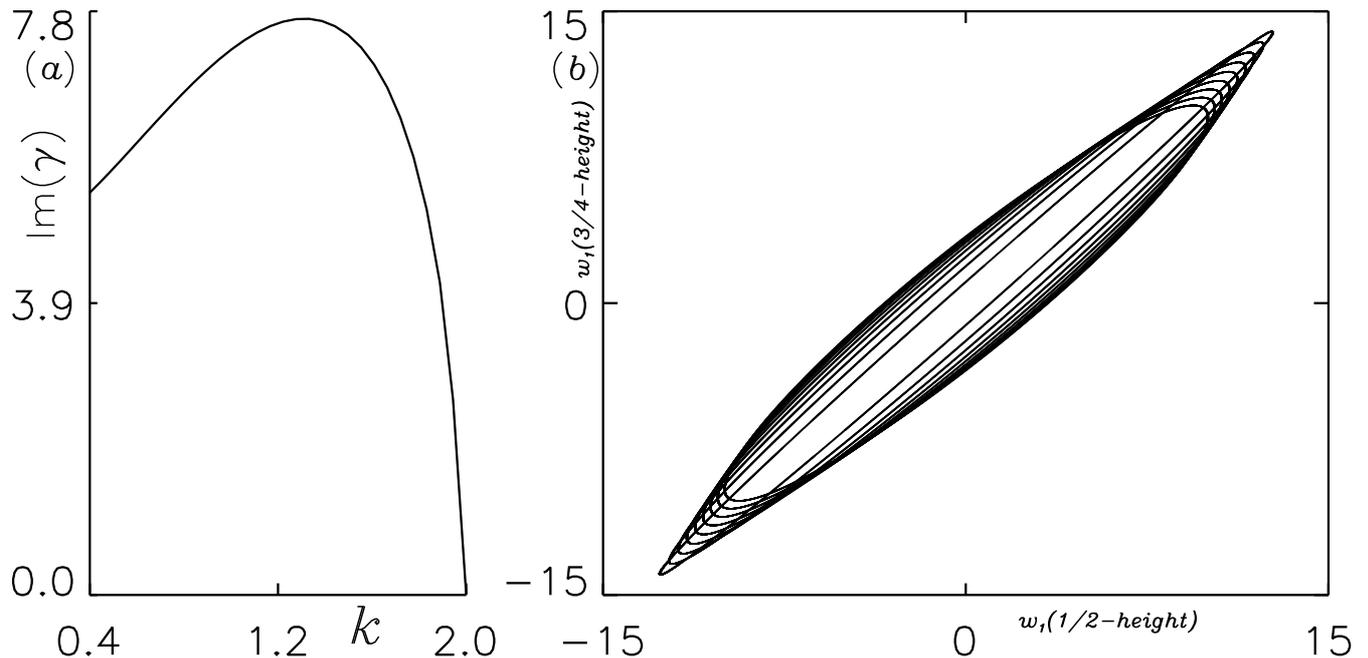}
\end{center}
\caption{(a) The positive frequency as a function of the wave number $k$ corresponding to the oscillatory branch of the neutral curve for $\theta_{i} =20$ shown in Fig. ~\ref{fig7} and (b) the corresponding phase portrait of perturbed fluid velocity component $w_1$. Fixed parameter values are $S_c = 20$, $V_c=15$, $\kappa=1$, $\mathrm{I_{D}}=0.5$, $\omega=0.4$, $\theta_{i}=20$ and 
${\mathrm{I}}_{t}=1$}
\label{portrait}
\end{figure}

\subsection*{Extinction coefficient $\kappa=1.0$} 
\label{os}

Here we consider the case, when $V_c = 15$, $\kappa=1$, $\mathrm{I_D}=0.5$ and $\omega=0.4$. At $\theta_{i}=0$, the location of the maximum basic concentration is around the mid-height of the chamber and the most unstable bioconvective solution remains in the stationary branch leading the solution to be stationary (non-oscillatory). As $\theta_{i}$ increases, the location of the maximum basic concentration shifts toward the top of the chamber. The location of the maximum basic concentration occurs at around $z = 0.65$ for $\theta_{i}=20$ [see Fig.~\ref{fig7}(a)] and  a single oscillatory branch bifurcates from the stationary branch at around wavenumber $k=2$ and defines a locus of points for $k \le 2.$ But, the most unstable solution still remains on the stationary branch, thus the solution is referred to stationary at bioconvective instability [see Fig.~\ref{fig7}(b)]. The qualitative nature of the oscillatory instability observed as above can be investigated by plotting the corresponding bifurcation diagram (phase portrait) because the bioconvective
motions become fully nonlinear on a timescale substantially less than the predicted period of oscillatory instability. 

\begin{figure}[!h]
\begin{center}
\includegraphics[scale=1]{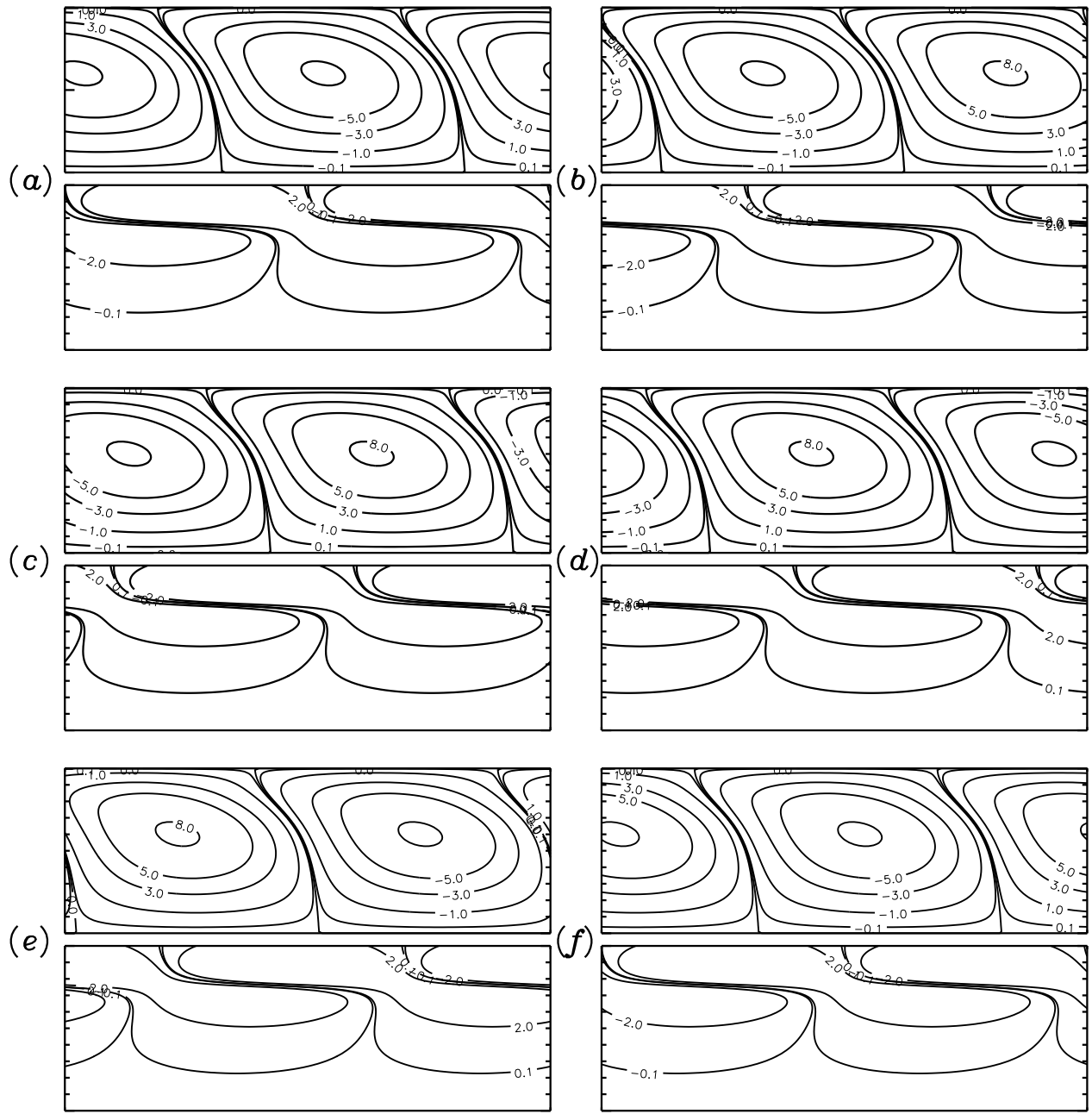}
\end{center}
\caption{Flow patterns of the perturbed velocity component $w_1$ during a cycle of oscillation at the onset of overstability.  The intervals between plots (a)-(f) are equal and the period is approximately $0.52$ units. Fixed parameter values are $S_c = 20$, $V_c=15$, $\kappa=1$, $\mathrm{I_D}=0.5$, $\omega=0.4$, $\theta_{i}=40$ and $\mathrm{I}_t=1$.}
\label{flow_pattern}
\end{figure} 

The qualitative behavior of the bioconvective system can be checked by plotting the corresponding phase portrait of infinitesimally small perturbations. Figure \ref{portrait}(a) shows the dependence of the positive frequency on the wavenumber $k$ on the oscillatory branch shown in Fig. ~\ref{fig7}(b) for  $\theta_{i}=20$. Since the eigenvalues $\gamma$ appear in complex conjugate pairs, only the positive frequency $\textnormal{Im}(\gamma)$ is shown in Fig. \ref{portrait}(a). When the frequency is decreased to zero, the oscillatory mode of disturbance changes to the stationary one at the onset of bioconvection [see Fig. \ref{portrait}(a)]. Since $\textnormal{Im}(\gamma) \neq 0$, the period of the oscillation, i.e. $2 \pi/\textnormal{Im}(\gamma),$  is the control (bifurcation) parameter. It varies by multiplying with an integer (i.e. $m 2 \pi /\textnormal{Im}(\gamma), \, m \in Z$) and the behaviour of the bioconvective system leads to a limit cycle or  periodic orbit. Fig. \ref{portrait}(b) shows that the bifurcation diagram has a degrading orbit with sequentially smaller radii upon each orbit, which implies that the bioconvective flow regime has damped oscillations as long as the frequency remains nonzero (i.e. when $k<2$)[see Fig. \ref{portrait}(a)]. The bioconvective flow regime experiences a steady convection as the frequency approaches to zero (i.e. when $k \ge 2$). 

\begin{figure}[!h]
\begin{center}
\includegraphics[scale=1]{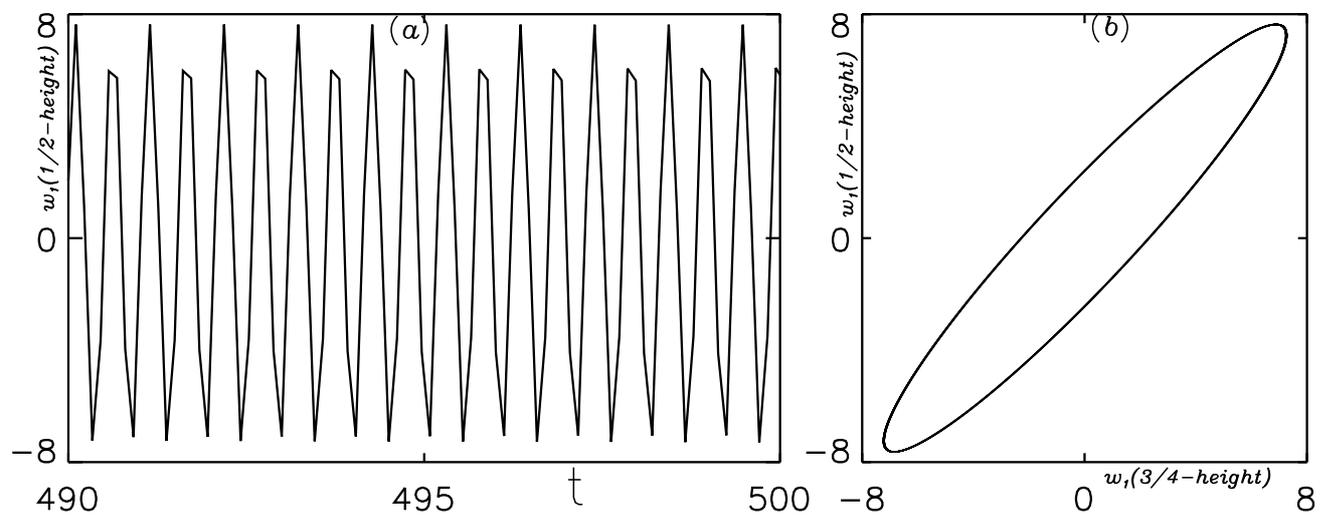}
\end{center}
\caption{(a) Time-evolving perturbed fluid velocity $w_1$ and (b) limit cycle. Fixed parameter values are $S_c = 20$, $V_c=15$, $\kappa=1$, $\mathrm{I_D}=0.5$, $\omega=0.4$, $\theta_{i}=40$ and $\mathrm{I}_t=1$.}
\label{limit_cycle}
\end{figure}

At $\theta_{i}=40$, the cells accumulate at around $z=0.8$ in the basic steady state [see Fig.~\ref{fig7}(a)]. In this instance, a single oscillatory branch bifurcates from the basic state at around $k_b \approx 3.4$ and the  oscillatory branch retains the most unstable solution making the bioconvective solution to be overstable [see Fig.~\ref{fig7}(b)]. Thus, the onset of overstability is at $k_c\approx 2.22$ and $R_c \approx 329.53$. At this point, two complex conjugate eigenvalues $\gamma =\pm 12.07i$ with zero real parts are found. This transition is known as Hopf bifurcation. The flow patterns of these two solutions corresponding to the conjugate pair of eigenvalues are mirror images of each other. The period of oscillation is $\dfrac{2 \pi}{\mathrm{Im}(\gamma)} \approx 0.52$ units. The bioconvective motions become fully nonlinear on a timescale substantially less than the predicted period of overstability. Hence, the flow patterns during one cycle of oscillation can be observed [see Fig. \ref{flow_pattern}]. This corresponds to a traveling wave solution moving toward the left of the figure. Fig. \ref{limit_cycle} shows the predicted time-varying perturbed fluid velocity component $w_1$ (Fig. \ref{limit_cycle}(a)) and its corresponding phase portrait (Fig. \ref{limit_cycle}(b)) at $k_c\approx 2.22$. Kindly note that here the period of the oscillation, $2\,\pi / \mathrm{Im}(\gamma)$ is the bifurcation parameter and hence the flow destabilization gives birth to a limit cycle [see Fig. \ref{limit_cycle}(b)]. The birth of a limit cycle due to the flow destabilization is again accepted as the Hopf bifurcation from bifurcation analysis. The supercritical nature (from the linear stability theory) of this Hopf bifurcation finally leads it into a stable limit cycle. The occurrence of overstability mode of disturbance in the bioconvective solution is again repeated when $\theta_{i}$ becomes $60$ [see Fig.~\ref{fig7}(b)]. An oscillatory branch still bifurcates from the stationary branch  when $\theta_{i}$ is increased to $80$ for which the location of the maximum concentration is at around the top of the chamber. But the unstable mode remains on the stationary branch [see Fig.~\ref{fig7}(b)]. Usually, the value of the critical Rayleigh number and wave number decreases as $\theta_{i}$ increases from zero for different fixed governing parameters. 

In this section, the bioconvective solution at instability for each value of $\theta_{i}$ is of mode $1$ except when $\theta_{i}$ is zero. The numerical results for critical Rayleigh number ($R_c$) and wavenumber ($k_c$) of this section are summarized in Table \ref{nr1}.

\begin{table}[!h]
\begin{center}
\caption{Bioconvective solutions with the variation of the diffuse irradiation by keeping other governing parameters fixed. A  result with double dagger symbol indicates that a smaller minimum occurs on an oscillatory branch and a starred result indicates that $R^{(1)}(k)$ branch of the neutral curve is oscillatory. \label{nr1}}
\begin{tabular}{c c c c c c c c c}
\hline \hline
$V_c$ \qquad& $\omega$  \qquad& $\kappa$ \qquad&  $\mathrm{I_D}$  \qquad&  $\theta_{i}(\textnormal{deg})$ \qquad&  $\lambda_c$  \qquad&   $R_c$   \qquad&   $\textnormal{Im}(\gamma)$   \qquad&  Mode \\[0.5pt]
$15$  \qquad&  $0.4$  \qquad& $0.5$  \qquad&  $0.26$  \qquad&  $0$ \qquad& $2.21$  \qquad&   $709.69$   \qquad&    $0$   \qquad&  $2$ \\[0.5pt]
$15$  \qquad&  $0.4$  \qquad& $0.5$  \qquad&  $0.26$  \qquad&  $20$  \qquad& $2.45$  \qquad&   $493.61$   \qquad&    $0$   \qquad&  $2$ \\[0.5pt]
$15$  \qquad&  $0.4$  \qquad& $0.5$  \qquad&  $0.26$  \qquad&  $40^{\star}$ \qquad& $2.82$  \qquad&   $266.36$   \qquad&    $0$   \qquad&  $1$ \\[0.5pt]
$15$  \qquad&  $0.4$  \qquad& $0.5$  \qquad&  $0.26$  \qquad&  $60$  \qquad& $2.56$  \qquad&   $232.97$   \qquad&    $0$   \qquad&  $1$ \\[0.5pt]
$15$  \qquad&  $0.4$  \qquad& $0.5$  \qquad&  $0.26$  \qquad&   $80$ \qquad& $2.45$  \qquad&   $242.59$   \qquad&    $0$   \qquad&  $1$ \\[0.5pt]
$15$  \qquad&  $0.4$  \qquad& $1.0$  \qquad&  $0.5$  \qquad&    $0$ \qquad& $1.8$  \qquad&   $668.55$   \qquad&    $0$   \qquad&  $2$ \\[0.5pt]
$15$  \qquad&  $0.4$  \qquad& $1.0$  \qquad&  $0.5$  \qquad&    $20^{\star}$ \qquad& $2.0$  \qquad&   $368.71$   \qquad&    $0$   \qquad&  $1$ \\[0.5pt]
$15$  \qquad&  $0.4$  \qquad& $1.0$  \qquad&  $0.5$   \qquad&   $40$ \qquad&   $2.82^\ddag$  \qquad&   $329.53^\ddag$   \qquad&    $12.07$   \qquad&  $1$ \\[0.5pt]
$15$  \qquad&  $0.4$  \qquad& $1.0$  \qquad&  $0.5$  \qquad&     $60$  \qquad& ${2.4}^\ddag$  \qquad&   $422.81^\ddag$   \qquad&  $13.84$   \qquad&  $1$ \\[0.5pt]
$15$  \qquad&  $0.4$  \qquad& $1.0$  \qquad&  $0.5$  \qquad&     $80^{\star}$ \qquad& $1.58$  \qquad&   $513$   \qquad&    $0$   \qquad&  $1$ \\[0.5pt]
\hline \hline
\end{tabular}
\end{center} 
\end{table}

\subsubsection{$V_c=10$ \textnormal{and}  $20$} 
\label{}

The bioconvective solutions due to the effects of oblique collimated irradiation at instability for $V_c=10$ and  $20$  are qualitatively similar to those of $V_c = 15$. The  numerical results for critical Rayleigh number ($R_c$) and wavelength ($\lambda_c$) of this section are presented in Table~\ref{nr2}.

\begin{table}[!h]
\begin{center}
\caption{Bioconvective solutions with the variation of the diffuse irradiation by keeping other governing parameters fixed. A  result with double dagger symbol indicates that a smaller minimum occurs on an oscillatory branch and a starred result indicates that $R^{(1)}(k)$ branch of the neutral curve is oscillatory. \label{nr2}}
\begin{tabular}{c c c c c c c c c}
\hline \hline
$V_c$ \qquad& $\omega$  \qquad& $\kappa$ \qquad&  $\mathrm{I_D}$  \qquad&  $\theta_{i}(\textnormal{deg})$ \qquad&  $\lambda_c$  \qquad&   $R_c$   \qquad&   $\textnormal{Im}(\gamma)$   \qquad&  Mode \\[0.5pt]
$10$  \qquad&  $0.4$  \qquad& $0.5$  \qquad&  $0.25$  \qquad&  $0$ \qquad& $2.45$  \qquad&   $946.08$   \qquad&    $0$   \qquad&  $2$ \\[0.5pt]
$10$  \qquad&  $0.4$  \qquad& $0.5$  \qquad&  $0.25$  \qquad&  $20$  \qquad& $2.75$  \qquad&   $724.75$   \qquad&    $0$   \qquad&  $1$ \\[0.5pt]
$10$  \qquad&  $0.4$  \qquad& $0.5$  \qquad&  $0.25$  \qquad&  $40$ \qquad& $3.88$  \qquad&   $356.06$   \qquad&    $0$   \qquad&  $1$ \\[0.5pt]
$10$  \qquad&  $0.4$  \qquad& $0.5$  \qquad&  $0.25$  \qquad&  $60$  \qquad& $4.17$  \qquad&   $216.86$   \qquad&    $0$   \qquad&  $1$ \\[0.5pt]
$10$  \qquad&  $0.4$  \qquad& $0.5$  \qquad&  $0.25$  \qquad&   $80$ \qquad& $3.88$  \qquad&   $189.45$   \qquad&    $0$   \qquad&  $1$ \\[0.5pt]
$10$  \qquad&  $0.4$  \qquad& $1.0$  \qquad&  $0.5$  \qquad&    $0$ \qquad& $1.97$  \qquad&   $851.87$   \qquad&    $0$   \qquad&  $2$ \\[0.5pt]
$10$  \qquad&  $0.4$  \qquad& $1.0$  \qquad&  $0.5$  \qquad&    $20$ \qquad& $2.16$  \qquad&   $588.25$   \qquad&    $0$   \qquad&  $1$ \\[0.5pt]
$10$  \qquad&  $0.4$  \qquad& $1.0$  \qquad&  $0.5$   \qquad&   $40^{\star}$ \qquad&   $2.3$  \qquad&   $344.25$   \qquad&    $0$   \qquad&  $1$ \\[0.5pt]
$10$  \qquad&  $0.4$  \qquad& $1.0$  \qquad&  $0.5$  \qquad&     $60^{\star}$  \qquad& ${2.16}$  \qquad&   $345$   \qquad&  $8.96$   \qquad&  $1$ \\[0.5pt]
$10$  \qquad&  $0.4$  \qquad& $1.0$  \qquad&  $0.5$  \qquad&     $80^{\star}$ \qquad& $2.16$  \qquad&   $308.33$   \qquad&    $0$   \qquad&  $1$ \\[0.5pt]
$20$  \qquad&  $0.4$  \qquad& $0.5$  \qquad&  $0.265$  \qquad&  $0$ \qquad& $2.0$  \qquad&   $659.34$   \qquad&    $0$   \qquad&  $2$ \\[0.5pt]
$20$  \qquad&  $0.4$  \qquad& $0.5$  \qquad&  $0.265$  \qquad&  $20^{\star}$  \qquad& $2.3$  \qquad&   $378.64$   \qquad&    $0$   \qquad&  $1$ \\[0.5pt]
$20$  \qquad&  $0.4$  \qquad& $0.5$  \qquad&  $0.265$  \qquad&  $40^{\star}$ \qquad& $2.19$  \qquad&   $276.8$   \qquad&    $0$   \qquad&  $1$ \\[0.5pt]
$20$  \qquad&  $0.4$  \qquad& $0.5$  \qquad&  $0.265$  \qquad&  $60^{\star}$  \qquad& $1.96$  \qquad&   $316.78$   \qquad&    $0$   \qquad&  $1$ \\[0.5pt]
$20$  \qquad&  $0.4$  \qquad& $0.5$  \qquad&  $0.265$  \qquad&   $80^{\star}$ \qquad& $1.86$  \qquad&   $359$   \qquad&    $0$   \qquad&  $1$ \\[0.5pt]
$20$  \qquad&  $0.4$  \qquad& $1.0$  \qquad&  $0.5$  \qquad&    $0$ \qquad& $1.53$  \qquad&   $826.83$   \qquad&    $0$   \qquad&  $2$ \\[0.5pt]
$20$  \qquad&  $0.4$  \qquad& $1.0$  \qquad&  $0.5$  \qquad&    $20^{\star}$ \qquad& $3.59$  \qquad&   $264.1$   \qquad&    $0$   \qquad&  $1$ \\[0.5pt]
$20$  \qquad&  $0.4$  \qquad& $1.0$  \qquad&  $0.5$   \qquad&   $40$ \qquad&   $2.62^\ddag$  \qquad&   $341.47^\ddag$   \qquad&    $22.22$   \qquad&  $1$ \\[0.5pt]
$20$  \qquad&  $0.4$  \qquad& $1.0$  \qquad&  $0.5$  \qquad&     $60$  \qquad& ${2.21}^\ddag$  \qquad&   $533.26^\ddag$   \qquad&  $27.14$   \qquad&  $1$ \\[0.5pt]
$20$  \qquad&  $0.4$  \qquad& $1.0$  \qquad&  $0.5$  \qquad&     $80$ \qquad& $2^\ddag$  \qquad&   $721.13^\ddag$   \qquad&    $23.44$   \qquad&  $1$ \\[0.5pt]
\hline \hline
\end{tabular}
\end{center} 
\end{table}

\subsubsection{Effects of swimming speed}

\begin{figure}[!h]
\begin{center}
\includegraphics[scale=1]{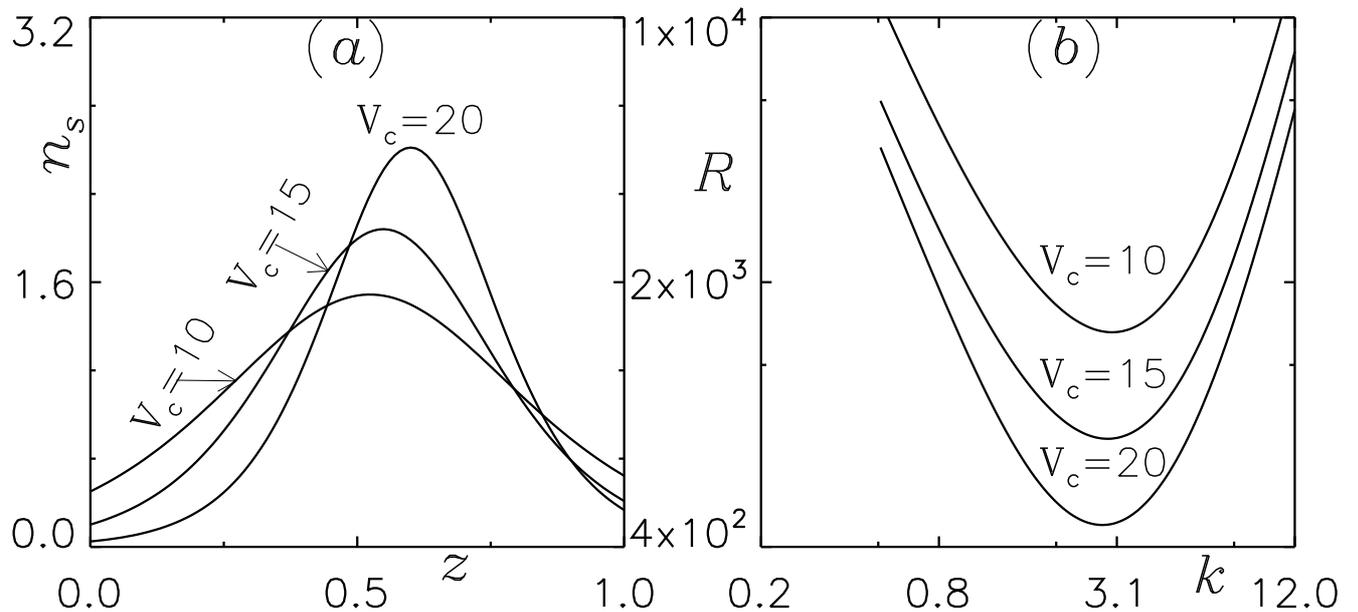}
\end{center}
\caption{(a) Base concentration profiles and(b) the corresponding neutral curves as swimming speed $V_c$ is increased.  Fixed parameter values are $S_c = 20$,  $\kappa=0.5$, $\mathrm{I_D}=0.26$, $\omega=0.4$, $\theta_{i}=0$, $G_{c}=1.3$ and $\mathrm{I}_t=1$.}
\label{eosp1}
\end{figure}

\begin{figure}[!h]
\begin{center}
\includegraphics[scale=1]{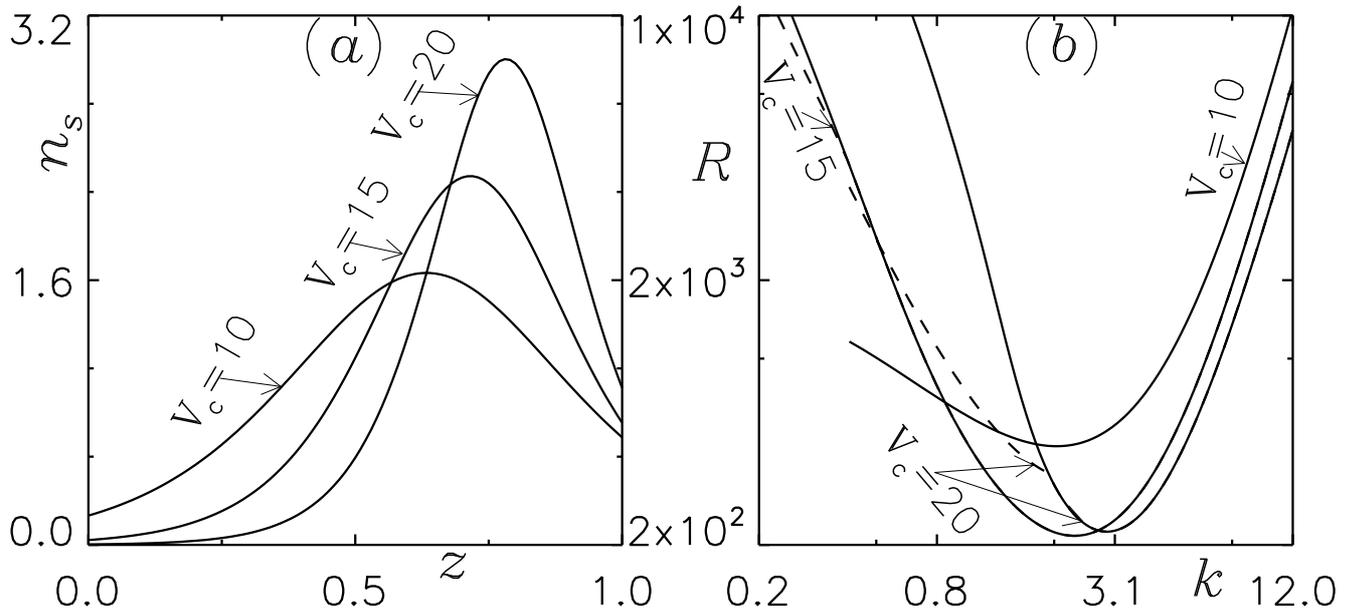}
\end{center}
\caption{(a) Base concentration profiles and(b) the corresponding neutral curves as swimming speed $V_c$ is increased.  Fixed parameter values are $S_c = 20$,  $\kappa=0.5$, $\mathrm{I_D}=0.26$, $\omega=0.4$, $\theta_{i}=40$, $G_{c}=1.3$ and $\mathrm{I}_t=1$.}
\label{xeosp}
\end{figure}

\begin{figure}[!h]
\begin{center}
\includegraphics[scale=1]{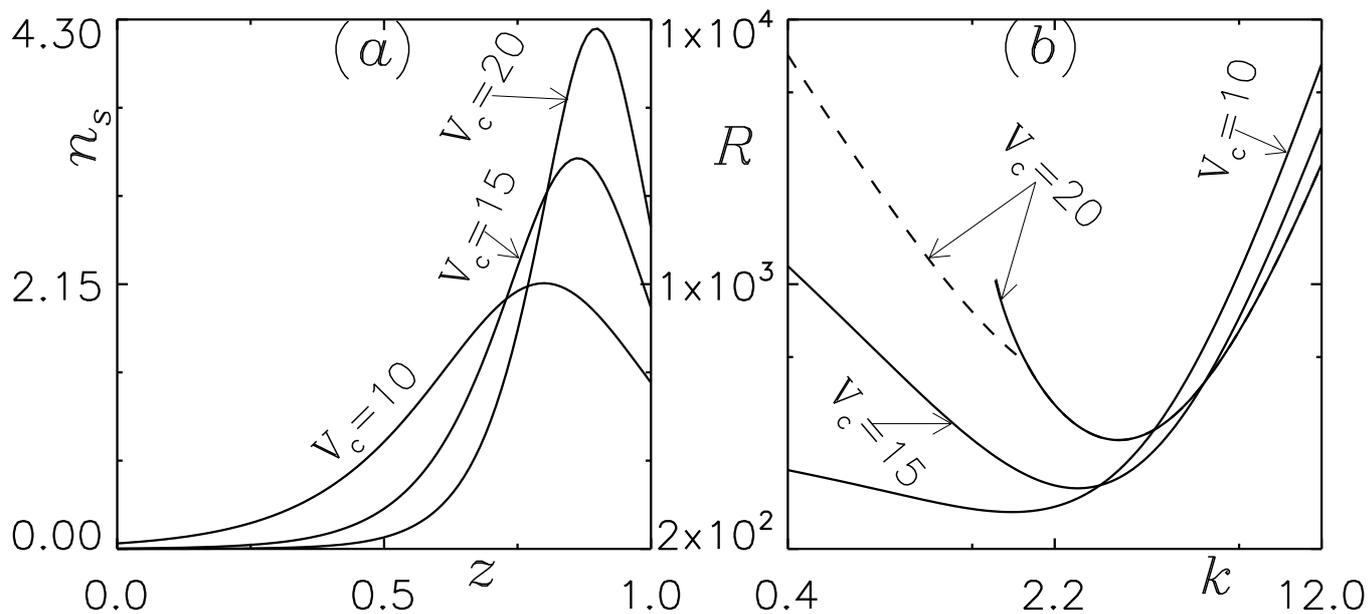}
\end{center}
\caption{(a) Base concentration profiles and(b) the corresponding neutral curves as swimming speed $V_c$ is increased.  Fixed parameter values are $S_c = 20$, $\kappa=0.5$, $\mathrm{I_D}=0.26$, $\omega=0.4$, $\theta_{i}=80$, $G_{c}=1.3$ and $\mathrm{I}_t=1$.}
\label{eosp2}
\end{figure}

Now we examine the effect of the cell swimming speed $V_c$ on the base concentration and the corresponding neutral curves at bioconvective instability by taking two different parameter ranges. We take a fixed set of parameters i.e. $\mathrm{I_D}=0.26$, $\omega=0.4$, $\kappa=0.5$, $\mathrm{I}_t=1$ and $\theta_{i}=0, 40, 80$. Next, we study the effect of $V_c$ on the solution of the linear stability problem by varying them discretely as $V_c=10, 15, 20$ and figures \ref{eosp1}, \ref{xeosp}, and \ref{eosp2} serve to illustrate them separately. The phototaxis function with critical intensity $G_c = 1.3$ [see Eq.~(\ref{taxis1})] is used in all the cases. Figure \ref{eosp1} shows the effect of $V_c$ on the critical states when the parameters  $\mathrm{I_D}=0.26$, $\omega=0.4$, $\kappa=0.5$, $\mathrm{I}_t=1$ and $\theta_{i}=0$ are kept fixed. In this case, the maximum base concentration occurs around the mid-height of the suspension when $V_{c}=10$. The location of maximum base concentration shifts towards the top when $V_c$ is increased to $15$ and $V_c=20$ respectively and the width (thickness) of the upper stable layer monotoncally decreases as $V_c$ is varied as $10$, $15$ and $20$ respectively. Thus, the effect of buoyancy of which tends to inhibit convective fluid motions decreases monotonically and simultaneously the critical wavenumber and critical Rayleigh number decrease as $V_c$ is increased as $10$, $15$ and $20$ respectively [see Fig. \ref{eosp1}].

Figure \ref{xeosp} shows the effect of the cell swimming speed $V_c$ on the basic concentration and the corresponding bioconvective instability (e.g critical wavelength and Rayleigh number) for fixed parameters $\kappa=0.5$, $\omega=0.4$, $\mathrm{I_D}=0.26$ and $\theta_{i}=40$.  In this case, the maximum concentration at the equilibrium state occurs at around $z \approx 0.6$ of the suspension for $V_c=10$ [see Fig.~\ref{xeosp}(a)]. As $V_c$ is increased  $15$,  the maximum concentration increases and its location occurs at $z \approx 0.7$. In this case, a single oscillatory branch bifurcates from the stationary branch at around $k \approx 0.4$ and defines a locus for $k \le 0.4.$ But, the stationary branch still retains the most unstable mode of disturbance. The critical Rayleigh number for $V_c=15$ is less as compared to the case when $V_c=10$ as the width of the stable layer decreases when $V_c=15$. Next, consider the case when $V_c=20.$ In this case, the maximum base concentration occurs around $z \approx 0.8$. In this instance, a single oscillatory branch bifurcates from the stationary branch at around $k \approx 1.88$ and defines a locus for $k \le 1.88$. But, the most unstable mode of disturbance remains on the stationary branch of the neutral curve. The thickness of the upper gravitationally stable region in the basic steady state is smallest for $V_c=20$ and it is largest for $V_c=10$ [see Fig.~\ref{xeosp}(a)]. Also, the steepness in base concentration gradient is largest for $V_c=20$ and it is smallest for $V_c=10$. The two factors, (more) steepness and (small) thickness of stable region, which support convection are more favourable for $V_c=20$ than at $V_c=10$. But, the positive phototaxis which opposes bioconvection is strongest for $V_c=20$ than at $V_c=10$. The latter feature dominates the former one and it results in an increase in critical Rayleigh number for $V_c=20$ than at $V_c=10$ [see Fig.~\ref{xeosp}(b)]. But, the former one dominates the later one for $V_c=15$ than $V_c=10$ and result in a lower critical Rayleigh number at $V_c=15$  than at $V_c=10$ [see Fig.~\ref{xeosp}(b)].

The effect of the cell swimming speed $V_c$ on the basic concentration profiles and the corresponding neutral curves for fixed parameters $\kappa=0.5$, $\omega=0.4$, $\mathrm{I_D}=0.26$ and $\theta_{i}=80$ are shown in Fig.~\ref{eosp2}. Here, the maximum base concentration for $V_c=10, 15, 20$ is located at around $z \approx 0.8$, $z \approx 0.85$, and $z \approx 0.9$  respectively of the suspension domain. The base concentration becomes steep for a higher swimming speed.  A steep concentration in the base steady state implies higher concentration gradient which supports the bioconvection, whereas the positive phototaxis offers higher resistance to the cells residing in the bioconvective plume at a higher swimming speed. Thus, the later effect dominates the higher gradient leading to higher critical Rayleigh number for a higher swimming speed [see Fig.~\ref{eosp2}(b)]. The perturbation to the basic steady state is stationary for $V_{c}=10, 15$. When $V_c=20,$ a single oscillatory branch bifurcates from the stationary branch at around $k \approx 1.85$ and defines a locus for $k \le 1.85$. But, the most unstable solution remains on the stationary branch of the neutral curve [see Fig.~\ref{eosp2}].

We have also investigated the effects of swimming speed on the bioconvective solutions when $\kappa=1$ and these effects remain qualitatively similar to the case of $\kappa=0.5$.  But, we observe oscillatory solutions for some parameters rather than stationary at bioconvective instability for a fixed $\theta{i}$ as $V_c$ increases from $10$ to $20$. Table~\ref{nr3} serves to illustrate all the results of this section.

\begin{table}[!h]
\begin{center}
\caption{Bioconvective solutions with the variation of the cell swimming speed by keeping  other governing parameters fixed. A  result with double dagger symbol indicates that a smaller minimum occurs on an oscillatory branch and a starred result indicates that $R^{(1)}(k)$ branch of the neutral curve is oscillatory. \label{nr3}}
\begin{tabular}{c c c c c c c c c}
\hline \hline
 $\mathrm{I_D}$ \qquad& $\omega$  \qquad& $\kappa$ \qquad&  $\theta_{i}(\textnormal{deg})$ \qquad&  $V_c$  \qquad&   $\lambda_c$  \qquad&   $R_c$   \qquad&   $\textnormal{Im}(\gamma)$   \qquad&  Mode \\[0.5pt]
$0.26$  \qquad&  $0.4$  \qquad& $0.5$  \qquad& $0$ \qquad& $10$  \qquad&   $2.11$  \qquad&   $1385.59$   \qquad&    $0$   \qquad&  $1$ \\[0.5pt]
$0.26$  \qquad&  $0.4$  \qquad& $0.5$  \qquad& $0$ \qquad& $15$  \qquad&   $2.19$  \qquad&   $709.68$   \qquad&    $0$   \qquad&  $1$ \\[0.5pt]
$0.26$  \qquad&  $0.4$  \qquad& $0.5$  \qquad& $0$ \qquad& $20$  \qquad&   $2.28$  \qquad&   $413.01$   \qquad&    $0$   \qquad&  $1$ \\[0.5pt]
$0.26$  \qquad&  $0.4$  \qquad& $0.5$  \qquad& $40$ \qquad& $10$  \qquad&   $3.23$  \qquad&   $497.17$   \qquad&    $0$   \qquad&  $1$ \\[0.5pt]
$0.26$  \qquad&  $0.4$  \qquad& $0.5$  \qquad& $40$ \qquad& $15^{\star}$  \qquad&   $2.8$  \qquad&   $266.33$   \qquad&    $0$   \qquad&  $1$ \\[0.5pt]
$0.26$  \qquad&  $0.4$  \qquad& $0.5$  \qquad& $40$ \qquad& $20^{\star}$  \qquad&   $2.16$  \qquad&   $274.14$   \qquad&    $0$   \qquad&  $1$ \\[0.5pt]
$0.26$  \qquad&  $0.4$  \qquad& $0.5$  \qquad& $80$ \qquad& $10$  \qquad&   $3.75$  \qquad&   $200.87$   \qquad&    $0$   \qquad&  $2$ \\[0.5pt]
$0.26$  \qquad&  $0.4$  \qquad& $0.5$  \qquad& $80$ \qquad& $15$  \qquad&   $2.49$  \qquad&   $242.56$   \qquad&    $0$   \qquad&  $2$ \\[0.5pt]
$0.26$  \qquad&  $0.4$  \qquad& $0.5$  \qquad& $80$ \qquad& $20^{\star}$  \qquad&   $1.89$  \qquad&   $355.65$   \qquad&    $0$   \qquad&  $2$ \\[0.5pt]
$0.48$  \qquad&  $0.4$  \qquad& $1.0$ \qquad& $0$ \qquad&  $10$  \qquad&   $2.25$  \qquad&   $523.41$   \qquad&    $0$   \qquad&  $1$ \\[0.5pt]
$0.48$  \qquad&  $0.4$  \qquad& $1.0$  \qquad& $0$ \qquad&  $15^{\star}$  \qquad&   $1.97$  \qquad&   $356.3$   \qquad&    $0$   \qquad&  $1$ \\[0.5pt]
$0.48$  \qquad&  $0.4$  \qquad& $1.0$  \qquad& $0$ \qquad&  $20$  \qquad&   $1.55$  \qquad&   $454.45$   \qquad&    $0$   \qquad&  $1$ \\[0.5pt]
$0.48$  \qquad&  $0.4$  \qquad& $1.0$  \qquad& $40$ \qquad&  $10^{\star}$  \qquad&   $2.28$  \qquad&   $316.53$   \qquad&    $0$   \qquad&  $1$ \\[0.5pt]
$0.48$  \qquad&  $0.4$  \qquad& $1.0$ \qquad& $40$ \qquad&  $15$  \qquad&   $2.68^{\ddag}$  \qquad&   $354.95^{\ddag}$   \qquad&    $12.59$   \qquad&  $1$ \\[0.5pt]
$0.48$  \qquad&  $0.4$  \qquad& $1.0$ \qquad& $40$ \qquad&  $20$  \qquad&   $2.4^{\ddag}$  \qquad&   $400.83^{\ddag}$   \qquad&    $23.98$   \qquad&  $1$ \\[0.5pt]
$0.48$  \qquad&  $0.4$  \qquad& $1.0$ \qquad& $80$ \qquad&  $10$  \qquad&   $2.25$  \qquad&   $297.46$   \qquad&    $0$   \qquad&  $2$ \\[0.5pt]
$0.48$  \qquad&  $0.4$  \qquad& $1.0$  \qquad& $80$ \qquad&  $15^{\star}$  \qquad&   $1.67$  \qquad&   $490.76$   \qquad&    $0$   \qquad&  $1$ \\[0.5pt]
$0.48$  \qquad&  $0.4$  \qquad& $1.0$ \qquad& $80$  \qquad&  $20$  \qquad&   $1.94^{\ddag}$  \qquad&   $793.03^{\ddag}$   \qquad&    $17.22$   \qquad&  $1$ \\[0.5pt]
\hline \hline
\end{tabular}
\end{center} 
\end{table}

\begin{figure}[!h]
\begin{center}
\includegraphics[scale=1]{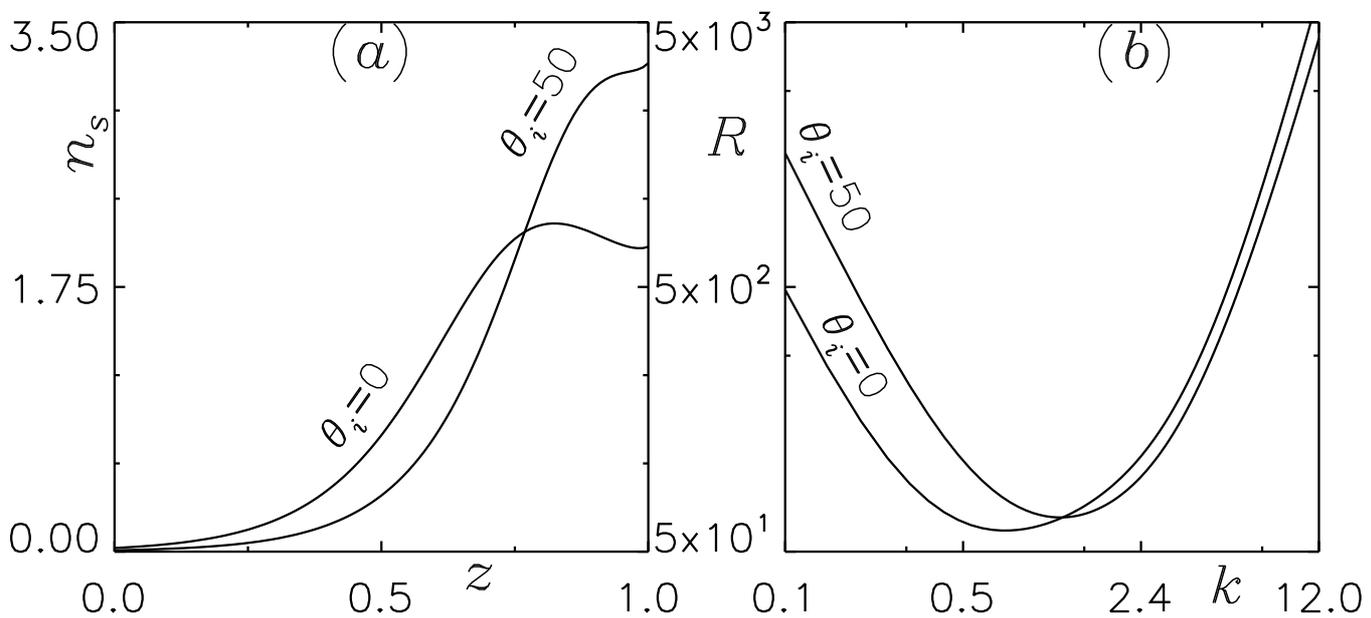}
\end{center}
\caption{(a) Base concentration profiles and (b) the corresponding neutral curves when $\theta_{i}=0,50$.  Fixed parameter values are $S_c = 20$, $\kappa=1$, $V_c=10$, $\omega=1$, $\mathrm{I_D}=0.02$, $G_{c}=1.9$ and ${\mathrm{I}}_{t}=1$.}
\label{ss01}
\end{figure}

\subsection{Case-II : Strong-scattering algal suspension} 
\label{di}

To emphasize the effects of oblique collimated irradiation on critical wavelength and Rayleigh number at bioconvective instability in a strong scattering algal suspension, here we consider the case when self-shading (absorption) is negligible by taking a purely scattering suspension i.e. $\omega=1$. We vary the parameter $V_c$ discretely as $V_c=10, 15, 20$ to study the effects of oblique collimated irradiation on bioconvection for fixed parameters $\kappa=1$, $\mathrm{I_D}=0.02$ and  $\omega=1$. In each case, we consider a phototaxis function with critical intensity $G_c=1.9$ whose mathematical form is given by  

%\begin{widetext}
\begin{equation}
T(G)=0.8\, \sin{\left(\dfrac{3\,\pi}{2} \Xi(G)\right)}-0.1\, \sin{\left(\dfrac{\pi}{2}\, \Xi(G)\right)}, \quad \Xi(G)=\dfrac{1}{3.8}\, G\, \exp{\left[0.135 \left(3.8-G\right)\right]}.
\label{taxis3}
\end{equation}
%\end{widetext}

\begin{figure}[!h]
\begin{center}
\includegraphics[scale=1]{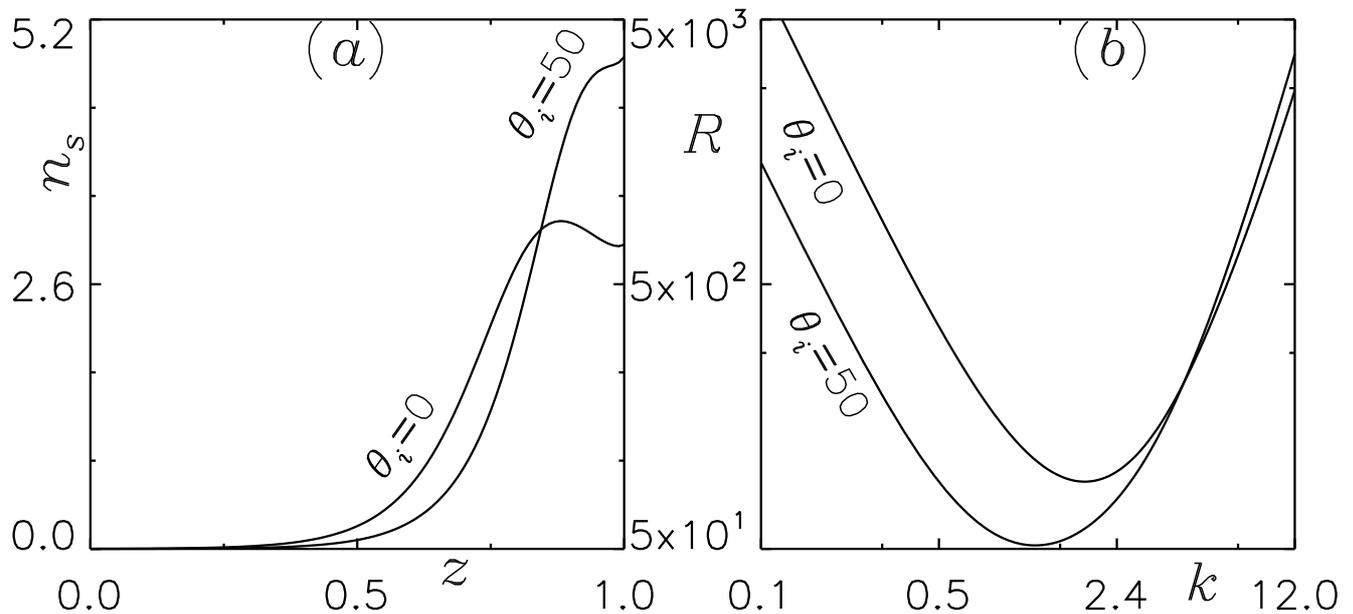}
\end{center}
\caption{(a) Base concentration profiles and (b) the corresponding neutral curves when $\theta_{i}=0,50$. Fixed parameter values are 
$S_c = 20$, $\kappa=1$, $V_c=15$, $\omega=1$, $\mathrm{I_D}=0.02$, $G_{c}=1.9$ and ${\mathrm{I}}_{t}=1$.}
\label{ss02}
\end{figure}

\begin{figure}[!h]
\begin{center}
\includegraphics[scale=1]{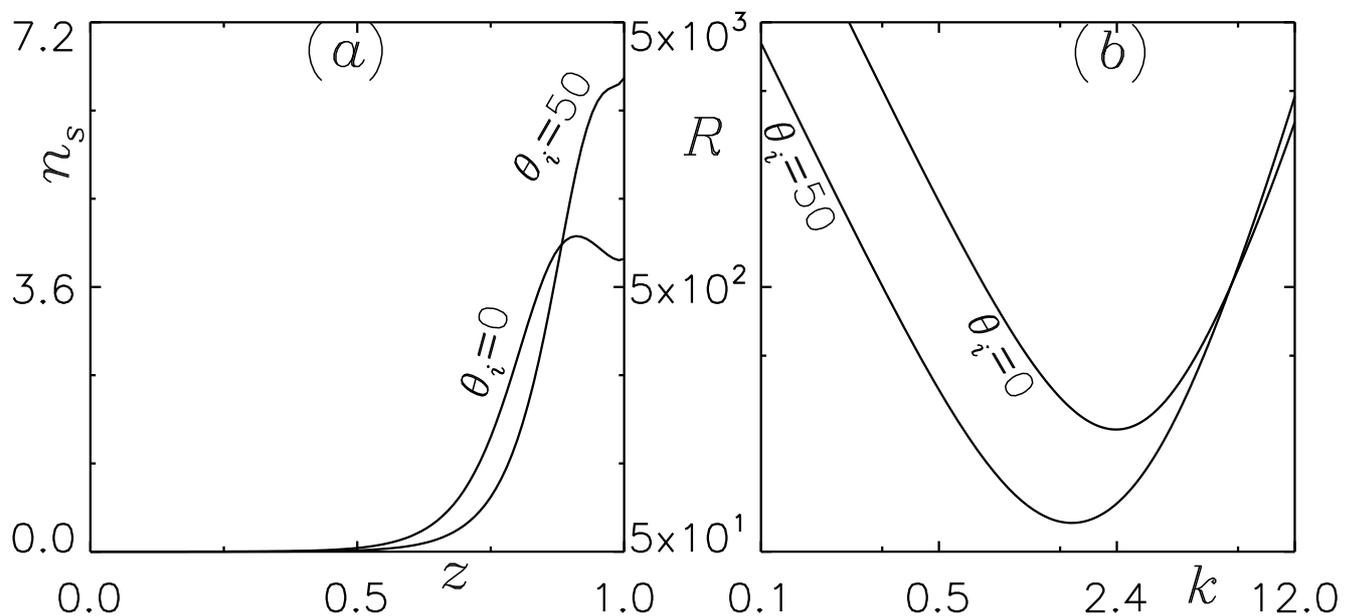}
\end{center}
\caption{(a) Base concentration profiles and (b) the corresponding neutral curves when $\theta_{i}=0, 50$.  Fixed parameter values are $S_c = 20$, $\kappa=1$, $V_c=20$, $\omega=1$, $\mathrm{I_D}=0.02$, $G_{c}=1.9$ and $\mathrm{I}_t=1$.}
\label{ss03}
\end{figure}

The unusual bimodal steady state is observed when $\omega=1$, $\theta_{i}=0$, $\kappa=1$ and $\mathrm{I_D}=0.02$ for the considered governing parameters in each case. The bimodal steady state in each case is transited to a unimodal equilibrium state when $\theta_{i}=50$. First, consider the case when $V_c=10$ for fixed $\omega=1$, $\theta_{i}=0$, $\kappa=1$ and $\mathrm{I_D}=0.02$ [see Fig.~\ref{ss01}]. Fig.~\ref{ss01}(a) and Fig.~\ref{ss01}(b) show the basic concentrations and the corresponding neutral stability curves for $\theta_{i}=0, 50$ respectively. In this case, the maximum base concentration occurs at two depths, $ z \approx 0.8$ and $z \approx 0.98$ of the suspension. Thus, negative (positive) phototaxis occurs inside (outside) the domain  between two locations $z \approx 0.8$ and $z \approx 0.98$ [see Fig.~\ref{ss01}(a)]. As $\theta_{i}$ increases to $50$, the negative phototaxis in the intermediate region between two depths, $ z \approx 0.8$ and $z \approx 0.98$ transits smoothly to positive phototaxis. Thus, the cells accumulate at around top of the suspension in the basic steady state [see Fig.~\ref{ss01}]. The base concentration profiles for the cases $V_c=15$ and $V_c=20$ are drawn in Fig. \ref{ss02}(a) and Fig. \ref{ss03}(a) and their behaviour remains similar to that seen when $V_c=10$ for $\theta_{i}=0, 50$. The width of the stable layer is same for $V_c=10, 15, 20$. The steepness in base concentration profile increases monotonically when $V_c$ varies from $10$ to $20$. The positive phototaxis offers higher resistance to the cells residing in the bioconvective plume at a higher swimming speed and the positive phototaxis which opposes bioconvection is stronger for $V_c =15$ than at $V_c=10$ and is strongest for $V_c = 20$ than at $V_c = 10,15$. But, a steep concentration in the base steady state implies higher concentration gradient which supports the bioconvection. The later effect dominates the former effect and thus the critical Rayleigh number decreases when $\theta_{i}$ is increased from $0$ to $50$ for $V_c=15, 20$ [see Fig. \ref{ss02} and Fig. \ref{ss03}].

\subsection{Comparison with up-swimming model}

Here we examine the effects of scattering on the linear stability of the isotropic scattering suspension illuminated by both diffuse and oblique collimated irradiation. We compare the observed results to those of up-swimming model of Panda \textit{et al.} \cite{PPS}. It is seen from Eq. (\ref{Tintensity}) that the basic equilibrium intensity is given by
%\begin{widetext}
\beq
G_s(z)=G_s^c(z)+G_s^d(z)=\mathrm{I}_t \, \exp\left(-\left(\dfrac{\kappa}{\cos{\theta_{0}}}\right)\int_z^1 n_s(z')\,dz'\right) +
G_s^d(z),  \label{ups}
\eeq
%\end{widetext}
where $G_s^d$ is the contribution due to isotropic scattering with diffuse and oblique collimated irradiation.

\begin{figure}[!h]
\begin{center}
\includegraphics[scale=1]{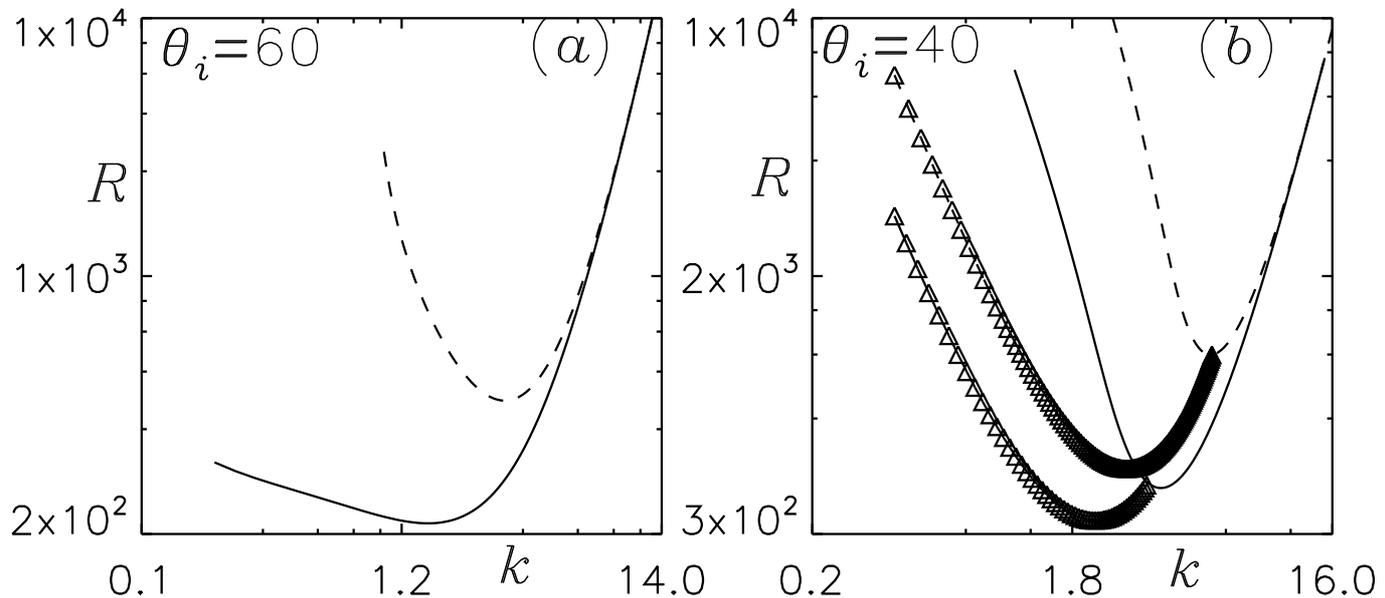}
\end{center}
\caption{Comparison of neutral curves for the model with scattering (solid lines) to that of up-swimming model (dashed lines). The overstable branches are marked with triangles. (a) corresponds to parameters $\theta_{i} = 60,$ $V_c=10$, $\kappa=0.5$, $\mathrm{I_D}=0.25$ and (b) corresponds to parameters $\theta_{i} = 40,$ $V_c=15$, $\kappa=1$,  $\mathrm{I_D}=0.5$ respectively.  Fixed parameter values are $S_c = 20,$ $G_{c}=1.3$, $\omega=0.4$ and $\mathrm{I}_t=1$.}
\label{fig8}
\end{figure}

Now we consider a suspension in which the intensity at a point $(x,y,z)$ is given by
\beq
G(x,y,z)=\mathrm{I}_t \, \exp\left(-\left(\dfrac{\kappa}{\cos{\theta_{0}}}\right)\int_z^1 n(x,y,z')\,dz'\right) + \chi(z),
\label{upsa}
\eeq
where $\chi(z)$ is independent of the concentration. If we take $\chi\equiv
G_s^d$, then the base flow will be the same but the cells swim in the vertical
direction only. Now the perturbed intensity  becomes
$$
G_1 = - \mathrm{I}_t \,\left(\dfrac{\kappa}{\cos{\theta_{0}}}\right) \left(\int_z^1n_1dz'\right)\exp\left(-\left(\dfrac{\kappa}{\cos{\theta_{0}}}\right)\int_z^1
n_s(z')\,dz'\right),  
$$
which is similar to that of Panda \textit{et al.} \cite{PPS}. Now using normal mode analsyis, the linear 
stability equation becomes
%\begin{widetext}
\beq
\left(\gamma S_c^{-1}+k^2-\frac{\ud^2}{\ud z^2}\right)\left(\frac{\ud^2}{\ud
z^2}-k^2\right)W=R k^2 D\Phi,\label{ulinwf}
\eeq     
\beq
\Gamma_1(z)\Phi +
\left[\gamma+k^2+\Gamma_2(z)\right]D\Phi +
V_c T_s\,D^2\Phi - D^3\Phi =-(Dn_s)W, \label{ulincf}  
\eeq  
%\end{widetext}
where
%\begin{widetext}
\begin{subequations}
\label{ugam}  
\begin{eqnarray}
\Gamma_1(z)&=& \left(\dfrac{\kappa}{\cos{\theta_{0}}}\right) V_c D\left(n_s G \frac{\ud T_s}{\ud G}\right), \label{ugb} \\
\Gamma_2(z)&=& 2\left(\dfrac{\kappa}{\cos{\theta_{0}}}\right) V_c n_s G_s^c \frac{\ud T_s}{\ud G}+V_c \frac{\ud T_s}{\ud G}D\chi+V_c \, n_s \, \left(\dfrac{\kappa}{\cos{\theta_{0}}}\right) \, \chi(z) \, \frac{\ud T_s}{\ud G}.
\label{ugc}
\end{eqnarray}
\end{subequations}  
%\end{widetext}
Thus, the basic flow is formed by the isotropic scattering suspension but the scattering is neglected in the perturbed equations. The phototaxis function with critical intensity $G_c = 1.3$ is used in these calculations [see Fig.~\ref{fig8}(a) and Eq.~(\ref{taxis1})] and the parameters $\omega=0.4$ and $\mathrm{I}_t=1$ are kept fixed.  In this section, we have included two sets of governing parameters based on stationary and overstability solutions at bioconvective instability.

We start with the stationary bioconvective solution for which the governing parameters $\theta_{i}=60$, $V_c=10$, $\kappa=0.5$ and $\mathrm{I_D}=0.25$ are kept fixed. In this case, the maximum base concentration is located around three-quarter height of the suspension and Figure~\ref{fig8}(a) depicts the comparison of neutral stability curves in between up-swimming \cite{PPS} and the present isotropic scattering model. Similarly,  Figure~\ref{fig8}(b) shows the comparison of overstability marginal  curves in between up-swimming \cite{PPS} and the present isotropic scattering model when the governing parameters  $\theta_{i}=40$, $V_c=15$, $\kappa=1$ and $\mathrm{I_D}=0.5$ are fixed. It reveals that the solutions of the linear stability problem from both the models agree well at small wavelegths and differ at large wavelengths too.

\section{CONCLUSIONS}
\label{chapt4:conclusions}

In this article, a model on phototaxis and self-shading that incorporates the effects of angle of incidence (or oblique collimated irradiation) on an isotropic scattering algal suspension is developed for the first time. The suspension is uniformly illuminated by both diffuse and oblique collimated irradiation. The onset of phototactic bioconvection via linear stability theory has been analyzed  using this model.

The obtained numerical results based on self-shading by the algae are summarized as follows. It is shown that the variation of total intensity across the suspension depth is not monotonic due to (isotropic) scattering by algae. Thus, at a higher scattering albedo (almost purely scattering suspension), the critical intensity locates at two different depths of the suspension . In this case, negative (positive) phototaxis occurs inside (outside) the intermediate region  between two locations where micro-organisms accumate in steady state. Thus, when $\theta_{i}$ is increased from zero to a certain nonzero value, the negative phototaxis which occurs in the region between two locations of the suspension converts to positive phototaxis. As a result, the bimodal steady state switches to a unimodal equilibrium state. It is worthy to note that an almost purely scattering suspension is very difficult to realize in quantitative studies on phototactic bioconvection as every suspension is absorbing to a large extent. 
 
The width of the upper stable layer decreases as the angle of incidence increases so that the effect of buoyancy which inhibits bioconvection monotonically decreases. Usually, the critical Wavenumber and critical Rayleigh number decrease  as the angle of incidence is increased from zero to higher non-zero values. The effect of the cell swimming speed on bioconvective solutions has been also investigated while keeping other governing parameters fixed. It is observed that the perturbation to the basic steady state transits from stationary to oscillatory state for a higher cell swimming speed.

In addition mode $2$ solutions occur on certain parts of the $R^{(1)}(k)$ branches of the neutral curve but the most unstable solution is always mode $1$ except, if the maximum basic concentration is located around mid-height of the suspension, then the most unstable solution can be mode $2$.

For certain ranges of parameter values, overstable bioconvective solutions are also observed, in which case the system undergoes a Hopf bifurcation at critical values, resulting in a travelling wave solution. Oscillatory bioconvective instability arises in many situations\cite{bst}, often when there is competition  between the stabilizing and destabilizing processes. We  see that three kinds of processes act in a phototactic suspension at the onset of bioconvective instability. The gravitationally stable region whose width decreases due to a non-zero higher value of $\theta_{i}$ above the layer of maximum concentration inhibits the bioconvection while the region below supports it. The role of phototaxis is twofold: it may inhibit or support the bioconvection.  The oscillatory/overstable bioconvective solutions may be observed due to these competing processes. 

The proposed model with (isotropic) scattering has been compared with the pure up-swimming model proposed by Panda \textit{et al.} \cite{PPS}. The models agree well at small wavelengths but differ at large wavelengths, because the contribution of scattering is negligible at small wavelengths. To study the effects of oblique collimated irradiation (light intensity)  on dominant initial pattern wavelength, we analyse the experimental results obtained by Williams and Bees \cite{williams_11}. It reveals that the initial wavelength of the instability increases with a decrease (an increase) in the light intensity (angle of incidence), which is in good agreement with the experimental results.

Quantitative study on bioconvection in a purely phototactic algal suspension is required for a comparison with the proposed work. But, unfortunately no such data exist till date to the best of our knowledge as most species of algae in a natural environment are gyrotactic or gravitactic in addition to being phototactic \cite{williams_11,hh:hh,khh:hhk}. The proposed model may be applied to other interesting problems involving populations of phototactic algae too. A suitable scattering phase function $p({\bmi s},{\bmi s}')$ would allow us to investigate the effects of oblique irradiation on a forward scattering algal suspension which is uniformly illuminated with/without diffuse irradiation.

\bibliographystyle{plain}

\end{document}